\documentclass[11pt]{article}

\usepackage[english]{babel} 
\usepackage[utf8]{inputenc} 
\usepackage[T1]{fontenc} 
\usepackage{fullpage}
\usepackage{amsfonts}
\usepackage{amsthm}
\usepackage{amsmath}
\usepackage{lmodern}
\usepackage{amssymb}
\usepackage{mathtools}
\usepackage{enumerate}
\usepackage{tikz-cd}
\usepackage{tikz}
\usepackage{graphicx}
\usepackage{stmaryrd}
\usepackage{longtable,graphics,multirow}
\usepackage{mathrsfs}
\usepackage{hyperref}
\usepackage{float}
\usepackage{colonequals}

\hypersetup{
colorlinks=true,
linkcolor=[rgb]{0,0.4,0.6}, %escala=255
filecolor=blue,
urlcolor=cyan,
citecolor=[rgb]{0,0.6,0},
anchorcolor=blue
}

\newtheorem{teor}{Theorem}
\numberwithin{teor}{section}

\newtheorem{lemma}[teor]{Lemma}

\newtheorem{prop}[teor]{Proposition}

\newtheorem{coro}[teor]{Corollary}

\newtheorem{remark}[teor]{Remark}

\theoremstyle{definition}\newtheorem{defi}[teor]{Definition}

\theoremstyle{definition}
\numberwithin{example}{subsection}

\newcommand{\fim}{ $_\blacksquare$}
\def\P{\mathbb{P}}

\def\Z{\mathbb{Z}}

\def\Q{\mathbb{Q}}

\begin{document}

\author{Renato Dias Costa}
\date{}
\title{Gaps on the intersection numbers of\\sections on a rational elliptic surface}
\maketitle

\begin{abstract}
%Let $X$ be a rational elliptic surface over an algebraically closed field with elliptic fibration $\pi:X\to\P^1$ and Mordell-Weil rank $r\geq 1$. Given $k\in\Bbb{N}$ we investigate if there is a pair of sections $P,P'$ with $P\cdot P'=k$ and say that $X$ has a $k$\textit{-gap} if not. We use Mordell-Weil lattices as our main tool and show how this is closely related to the problem of representing integers by quadratic forms. We describe all gaps on $X$ having torsion-free Mordell-Weil group and $r=1$. We prove that $X$ has no gaps if $r\geq 5$, while gaps occur with probability $1$ if $r=1,2$. We also show that $X$ has a $1$-gap if and only if $\pi$ has $\text{III}^*$ fiber and a $2$-gap if and only if $\pi$ has a $\text{I}^*_3$ fiber.
%\\ \\
\indent Given a rational elliptic surface $X$ over an algebraically closed field, we investigate whether a given natural number $k$ can be the intersection number of two sections of $X$. If not, we say that $k$ a \textit{gap number}. We try to answer when gap numbers exist, how they are distributed and how to identify them. We use Mordell-Weil lattices as our main tool, which connects the investigation to the classical problem of representing integers by positive-definite quadratic forms. 

%We explore conditions under which  numbers exist and, when they do, we try to describe their distribution and in some cases identify them. We use Mordell-Weil lattices as our main tool, which reveals a connection with the classical problem of representing integers by positive-definite quadratic forms.

%We explore necessary and sufficient conditions for finding a pair of sections with intersection number k and apply these results to answer questions about the existence of gap numbers, how they are distributed and, in some cases, how to identify them precisely identification.  
\end{abstract}

\tableofcontents

\section{Introduction}\label{section:intro}
\noindent\textbf{Description of the problem.} Let $X$ be a rational elliptic surface over an algebraically closed field, i.e. a smooth, rational projective surface with a fibration $\pi:X\to\P^1$ whose general fiber is a smooth curve of genus $1$. Assume also that $\pi$ is relatively minimal, i.e. no fiber contains an exceptional curve in its support. We use $E/K$ to denote the generic fiber of $\pi$, which is an elliptic curve over the function field $K:=k(\P^1)$. By the Mordell-Weil theorem, the set $E(K)$ of $K$-points is a finitely generated Abelian group, whose rank we denote by $r$. The points on $E(K)$ are in bijective correspondence with the sections of $\pi$, as well as with the exceptional curves on $X$, so we use these terms interchangeably. This paper addresses the following question: given sections $P_1,P_2\in E(K)$, what values can the intersection number $P_1\cdot P_2$ possibly attain?
\\ \\
\noindent \textbf{Original motivation.} The problem originates from a previous investigation of conic bundles on $X$, i.e. morphisms $\varphi:X\to\P^1$ whose general fiber is a smooth curve of genus zero \cite{Costa}. More specifically, one of the ways to produce a conic bundle is by finding a pair of sections $P_1,P_2\in E(K)$ with $P_1\cdot P_2=1$, so that the linear system $|P_1+P_2|$ induces a conic bundle $\varphi_{|P_1+P_2|}:X\to\P^1$ having $P_1+P_2$ as a reducible fiber. We may ask under which conditions such a pair exists. An immediate necessary condition is that $r\geq 1$, for if $r=0$ any two distinct sections must be disjoint \cite[Cor. 8.30]{MWL}. Conversely, given that $r\geq 1$, does $X$ admit such a pair? The first observation is that $r\geq 1$ implies an infinite number of sections, so we should expect infinitely many values for $P_1\cdot P_2$ as $P_1,P_2$ run through $E(K)$. Then the question is ultimately: what values may $P_1\cdot P_2$ assume?
\\ \\
\noindent \textbf{Mordell-Weil lattices.} The computation of intersection numbers on a surface is a difficult problem in general. However, as we are concerned with sections on an elliptic surface, the information we need is considerably more accessible. The reason for this lies in the Mordell-Weil lattice, a concept first established in \cite{Elkies}, \cite{Shioda}, \cite{Shioda'}. It involves the definition of a $\Q$-valued pairing $\langle\cdot,\cdot\rangle$ on $E(K)$, called the \textit{height pairing} \cite[Section 6.5]{MWL}, inducing a positive-definite lattice $(E(K)/E(K)_\text{tor},\langle\cdot,\cdot\rangle)$, named the Mordell-Weil lattice. A key aspect of its construction is the connection with the Néron-Severi lattice, so that the height pairing and the intersection pairing of sections are strongly intertwined. In the case of rational elliptic surfaces, the possibilities for the Mordell-Weil lattice have already been classified in \cite{OguisoShioda}, which gives us a good starting point.
\ \\
\noindent \textbf{Representation of integers.} The use of Mordell-Weil lattices in our investigation leads to a classical problem in number theory, which is the representation of integers by positive-definite quadratic forms. Indeed, the free part of $E(K)$ is generated by $r$ terms, so the height $h(P):=\langle P,P\rangle$ induces a positive-definite quadratic form on $r$ variables with coefficients in $\Bbb{Q}$. If $O\in E(K)$ is the neutral section and $R$ is the set of reducible fibers of $\pi$, then by the height formula (\ref{height_formula_P})
$$h(P)=2+2(P\cdot O)-\sum_{v\in R}\text{contr}_v(P),$$
where the sum over $v$ is a rational number which can be estimated. By clearing denominators, we see that the possible values of $P\cdot O$ depend on a certain range of integers represented by a positive-definite quadratic form with coefficients in $\Z$. This point of view is explored in some parts of the paper, where we apply results such as the classical Lagrange's four-square theorem \cite[\S 20.5]{HardyWright}, the counting of integers represented by a binary quadratic form \cite[p. 91]{Bernays} and the more recent Bhargava-Hanke's 290-theorem on universal quadratic forms \cite[Thm. 1]{BhargavaHanke}.
\\ \\
\textbf{Statement of results.} Given $k\in\Bbb{Z}_{\geq 0}$ we investigate whether there is a pair of sections $P_1,P_2\in E(K)$ such that $P_1\cdot P_2=k$. If such a pair does not exist, we say that $X$ has a $k$\textit{-gap}, or that $k$ is a \textit{gap number}. Our first result is a complete identification of gap numbers in some cases:
\\ \\
\noindent \textbf{Theorem \ref{thm:identification_of_gaps_r=1}.} If $E(K)$ is torsion-free with rank $r=1$, we have the following characterization of gap numbers on $X$ according to the lattice $T$ associated to the reducible fibers of $\pi$.
\begin{table}[h]
\begin{center}
\centering
\begin{tabular}{cc} 
\hline
\multirow{2}{*}{$T$} & \multirow{2}{*}{$\begin{matrix}k\text{ is a gap number}\Leftrightarrow\text{none of}\\ \text{the following are perfect squares}\end{matrix}$}\\ 
& \\
\hline
\multirow{2}{*}{$E_7$} & \multirow{2}{*}{$k+1$, $4k+1$}\\ 
& \\
\hline
\multirow{2}{*}{$A_7$} & \multirow{2}{*}{$\frac{k+1}{4}$, $16k,...,16k+9$}\\ 
& \\
\hline
\multirow{2}{*}{$D_7$} & \multirow{2}{*}{$\frac{k+1}{2}$, $8k+1,...,8k+4$}\\ 
& \\
\hline
\multirow{2}{*}{$A_6\oplus A_1$} & \multirow{2}{*}{$\frac{k+1}{7}$, $28k-3,...,28k+21$}\\ 
& \\
\hline
\multirow{2}{*}{$E_6\oplus A_1$} & \multirow{2}{*}{$\frac{k+1}{3}$, $12k+1,...,12k+9$}\\ 
& \\
\hline
\multirow{2}{*}{$D_5\oplus A_2$} & \multirow{2}{*}{$\frac{k+1}{6}$, $24k+1,...,24k+16$}\\ 
& \\
\hline
\multirow{2}{*}{$A_4\oplus A_3$} & \multirow{2}{*}{$\frac{k+1}{10}$, $40k-4,...,40k+25$}\\ 
& \\
\hline
\multirow{2}{*}{$A_4\oplus A_2\oplus A_1$} & \multirow{2}{*}{$\frac{k+1}{15}$, $60k-11,...,60k+45$}\\ 
& \\
\hline
\end{tabular}
\end{center}
\end{table}

We also explore the possibility of $X$ having no gap numbers. We prove that, in fact, this is always the case if the Mordell-Weil rank is big enough.
%. Applying the classical Lagrange's four-square theorem \cite[Thm. 10.6]{JJ} and the more recent Bhargava-Hanke's 290-theorem on universal quadratic forms \cite[Thm. 1]{BhargavaHanke}, we show that $X$ has no gaps if the Mordell-Weil rank is big enough.
\\ \\
\noindent\textbf{Theorem \ref{thm:gap_free_r_geq_5}.} If $r\geq 5$, then $X$ has no gap numbers.
\\ \\
\indent On the other hand, for $r\leq 2$ we show that gap numbers occur with probability $1$.
\\ \\
\noindent\textbf{Theorem \ref{thm:gaps_probability_1_r=1,2}.} If $r\leq 2$, then the set of gap numbers of $X$, i.e. $G:=\{k\in\Bbb{N}\mid k\text{ is a gap number of }X\}$ has density $1$ in $\Bbb{N}$, i.e.
$$\lim_{n\to\infty}\frac{\#G\cap\{1,...,n\}}{n}=1.$$

\indent At last we answer the question from the original motivation, which consists in classifying the rational elliptic surfaces with a $1$-gap:
\\ \\
\noindent\textbf{Theorem \ref{thm:surfaces_with_a_1-gap}.} $X$ has a $1$-gap if and only if $r=0$ or $r=1$ and $\pi$ has a $\text{III}^*$ fiber.
\newpage
%\\ \\
%\indent At last, using the same techniques we are able to classify the rational elliptic surfaces with a $2$-gap.
%\\ \\
%\noindent\textbf{Theorem \ref{gap_2}.} $X$ has a $2$-gap if and only if $r=0$ or $r=1$ and $\pi$ has a $\text{I}_3^*$ fiber.
%\\ \\
\noindent\textbf{Structure of the paper.} The text is organized as follows. Section~\ref{section:preliminaries} introduces the main objects, namely the Mordell-Weil lattice, the bounds $c_\text{max},c_\text{min}$ for the contribution term, the difference $\Delta=c_\text{max}-c_\text{min}$ and the quadratic form $Q_X$ induced by the height pairing. In Section~\ref{section:the_role_of_torsion_sections} we explain the role of torsion sections in the investigation. The key technical results are gathered in Section~\ref{section:existence_of_a_pair}, where we state necessary and sufficient conditions for having $P_1\cdot P_2=k$ for a given $k$. Section~\ref{section:main_results} contains the main results of the paper, namely: the description of gap numbers when $E(K)$ is torsion-free with $r=1$ (Subsection~\ref{subsection:identification_of_gaps_r=1}), the absence of gap numbers for $r\geq 5$ (Subsection~\ref{subsection:gap_free_r_geq_5}), density of gap numbers when $r\leq 2$ (Subsection~\ref{subsection:gaps_probability_1_r=1,2}) and the classification of surfaces with a $1$-gap (Subsection~\ref{subsection:surfaces_with_a_1-gap}). Section~\ref{section:appendix} is an appendix containing Table~\ref{tabela_completa}, which stores the relevant information about the Mordell-Weil lattices of rational elliptic surfaces with $r\geq 1$.

\section{Preliminaries}\label{section:preliminaries}\
\indent Throughout the paper $X$ denotes a rational elliptic surface over an algebraically closed field $k$ of any characteristic. More precisely, $X$ is a smooth rational projective surface with a fibration $\pi:X\to\P^1$, with a section, whose general fiber is a smooth curve of genus $1$. We assume moreover that $\pi$ is relatively minimal (i.e. each fiber has no exceptional curve in its support) \cite[Def. 5.2]{MWL}. The generic fiber of $\pi$ is an elliptic curve $E/K$ over $K:=k(\P^1)$. The set $E(K)$ of $K$-points is called the \textit{Mordell-Weil group} of $X$, whose rank is called the \textit{Mordell-Weil rank of }$X$, denoted by
$$r:=\text{rank }E(K).$$

In what follows we introduce the main objects of our investigation and stablish some notation. 

\subsection{The Mordell-Weil Lattice}\label{subsection:MW_lattice}\
\indent We give a brief description of the Mordell-Weil lattice, which is the central tool used in the paper. Although it can be defined on elliptic surfaces in general, we restrict ourselves to rational elliptic surfaces. For more information on Mordell-Weil lattices, we refer the reader to the comprehensive introduction by Schuett and Shioda \cite{MWL} in addition to the original sources, namely \cite{Elkies}, \cite{Shioda}, \cite{Shioda'}.

We begin by noting that points in $E(K)$ can be regarded as curves on $X$ and by defining the lattice $T$ and the trivial lattice $\text{Triv}(X)$, which are needed to define the Mordell-Weil lattice.
\\ \\
\noindent\textbf{Sections, points on $E(K)$ and exceptional curves.} The sections of $\pi$ are in bijective correspondence with points on $E(K)$. Moreover, since $X$ is rational and relatively minimal, points on $E(K)$ also correspond to exceptional curves on $X$ \cite[Section 8.2]{Schuett-Shioda}. For this reason we identify sections of $\pi$, points on $E(K)$ and exceptional curves on $X$. 
\\ \\
\noindent\textbf{The lattice $T$ and the trivial lattice $\text{Triv}(X)$.} Let $O\in E(K)$ be the neutral section and $R:=\{v\in\P^1\mid \pi^{-1}(v)\text{ is reducible}\}$ the set of reducible fibers of $\pi$. The components of a fiber $\pi^{-1}(v)$ are denoted by $\Theta_{v,i}$, where $\Theta_{v,0}$ is the only component intersected by $O$. The Néron-Severi group $\text{NS}(X)$ together with the intersection pairing is called the Néron-Severi lattice. 
\newpage
\indent We define the following sublattices of $\text{NS}(X)$, which encode the reducible fibers of $\pi$:
$$T_v:=\Bbb{Z}\langle \Theta_{v,i}\mid i\neq 0\rangle\text{ for }v\in R,$$
$$T:=\bigoplus_{v\in R}T_v.$$

By Kodaira's classification \cite[Thm. 5.12]{MWL}, each $T_v$ with $v\in R$ is represented by a Dynkin diagram $A_m,D_m$ or $E_m$ for some $m$. We also define the \textit{trivial lattice} of $X$, namely
$$\text{Triv}(X):=\Bbb{Z}\langle O, \Theta_{v,i}\mid i\geq 0,\,v\in R \rangle.$$

Next we define the Mordell-Weil lattice and present the height formula.
\\ \\
\noindent\textbf{The Mordell-Weil lattice.} In order to give $E(K)$ a lattice structure, we cannot use the intersection pairing directly, which only defines a lattice on $\text{NS}(X)$ but not on $E(K)$. This is achieved by defining a $\Q$-valued pairing, called the \textit{height pairing}, given by
\begin{align*}
\langle\cdot,\cdot\rangle:E(K)\times E(K)&\to\Q\\
P,Q&\mapsto -\varphi(P)\cdot\varphi(Q),
\end{align*}
where $\varphi:E(K)\to\text{NS}(X)\otimes_\Z \Q$ is defined from the orthogonal projection with respect to $\text{Triv}(X)$ (for a detailed exposition, see \cite[Section 6.5]{MWL}). Moreover, dividing by torsion elements we get a positive-definite lattice $(E(K)/E(K)_\text{tor},\langle \cdot,\cdot\rangle)$ \cite[Thm. 6.20]{MWL}, called the \textit{Mordell-Weil lattice}. 
\\ \\
\noindent\textbf{The height formula.} The height pairing can be explicitly computed by the \text{height formula} \cite[Thm. 6.24]{MWL}. For rational elliptic surfaces, it is given by
\begin{equation}\label{height_formula_PQ}
\langle P,Q\rangle=1+(P\cdot O)+(Q\cdot O)-(P\cdot Q)-\sum_{v\in R}\text{contr}_v(P,Q),
\end{equation}
\begin{equation}\label{height_formula_P}
h(P):=\langle P,P\rangle=2+2(P\cdot O)-\sum_{v\in R}\text{contr}_v(P),
\end{equation}

\noindent where $\text{contr}_v(P):=\text{contr}_v(P,P)$ and $\text{contr}_v(P,Q)$ are given by Table~\ref{contributions} \cite[Table 6.1]{MWL} assuming $P,Q$ meet $\pi^{-1}(v)$ at $\Theta_{v,i},\Theta_{v,j}$ resp. with $0<i<j$. If $P$ or $Q$ meets $\Theta_{v,0}$, then $\text{contr}_v(P,Q):=0$.

\begin{table}[h]
\begin{center}
\centering
\begin{tabular}{c|c|c|c|c|c|c} 
\hline
\multirow{2}{*}{\hfil $T_v$} & \multirow{2}{*}{\hfil $A_1$} & \multirow{2}{*}{\hfil $E_7$} & \multirow{2}{*}{\hfil $A_2$} & \multirow{2}{*}{\hfil $E_6$} & \multirow{2}{*}{\hfil $A_{n-1}$} & \multirow{2}{*}{\hfil $D_{n+4}$}\\ 
& & & & & &\\
\hline
\multirow{2}{*}{\hfil $\text{Type of }\pi^{-1}(v)$} & \multirow{2}{*}{\hfil $\text{III}$} & \multirow{2}{*}{\hfil $\text{III}^*$} & \multirow{2}{*}{\hfil $\text{IV}$} & \multirow{2}{*}{\hfil $\text{IV}^*$} & \multirow{2}{*}{\hfil $\text{I}_n$} & \multirow{2}{*}{\hfil $\text{I}_n^*$}\\ 
& & & & & &\\
\hline
\multirow{3}{*}{\hfil $\text{contr}_v(P)$} & \multirow{3}{*}{\hfil $\frac{1}{2}$} & \multirow{3}{*}{\hfil $\frac{3}{2}$} & \multirow{3}{*}{\hfil $\frac{2}{3}$} & \multirow{3}{*}{\hfil $\frac{4}{3}$} & \multirow{3}{*}{\hfil $\frac{i(n-i)}{n}$} & 	\multirow{3}{*}{\hfil $\begin{cases}1 & (i=1)\\1+\frac{n}{4} & (i>1)\end{cases}$}\\
& & & & & &\\
& & & & & &\\
\hline
\multirow{3}{*}{\hfil $\text{contr}_v(P,Q)$} & \multirow{3}{*}{\hfil $\text{-}$} & \multirow{3}{*}{\hfil $\text{-}$} & \multirow{3}{*}{\hfil $\frac{1}{3}$} & \multirow{3}{*}{\hfil $\frac{2}{3}$} & \multirow{3}{*}{\hfil $\frac{i(n-j)}{n}$} & \multirow{3}{*}{\hfil $\begin{cases}\frac{1}{2} & (i=1)\\\frac{1}{2}+\frac{n}{4} & (i>1)\end{cases}$}\\
& & & & & &\\
& & & & & &\\
\hline
\end{tabular}
\caption{Local contributions from reducible fibers to the height pairing.}\label{contributions}
\end{center}
\end{table}
\ \\
\noindent\textbf{The minimal norm.} Since $E(K)$ is finitely generated, there is a minimal positive value for $h(P)$ as $P$ runs through $E(K)$ with $h(P)>0$. It is called the \textit{minimal norm}, denoted by
$$\mu:=\min\{h(P)>0\mid P\in E(K)\}.$$

\noindent\textbf{The narrow Mordell-Weil lattice.} An important sublattice of $E(K)$ is the \textit{narrow Mordell-Weil lattice} $E(K)^0$, defined as
\begin{align*}
E(K)^0&:=\{P\in E(K)\mid P\text{ intersects }\Theta_{v,0}\text{ for all }v\in R\}\\
&=\{P\in E(K)\mid \text{contr}_v(P)=0\text{ for all } v\in R\}.
\end{align*}

As a subgroup, $E(K)^0$ is torsion-free; as a sublattice, it is a positive-definite even integral lattice with finite index in $E(K)$ \cite[Thm. 6.44]{MWL}. The importance of the narrow lattice can be explained by its considerable size as a sublattice and by the easiness to compute the height pairing on it, since all contribution terms vanish. A complete classification of the lattices $E(K)$ and $E(K)^0$ on rational elliptic surfaces is found in \cite[Main Thm.]{OguisoShioda}.

\subsection{Gap numbers}\label{subsection:gap_numbers}\
\indent We introduce some convenient terminology to express the possibility of finding a pair of sections with a given intersection number.

\begin{defi}
If there are no sections $P_1,P_2\in E(K)$ such that $P_1\cdot P_2=k$, we say that $X$ has a $k$-\textit{gap} or that $k$ is a \textit{gap number} of $X$.
\end{defi}

\begin{defi}
We say that $X$ is \textit{gap-free} if for every $k\in\Bbb{Z}_{\geq 0}$ there are sections $P_1,P_2\in E(K)$ such that $P_1\cdot P_2=k$.
\end{defi}

\begin{remark}\label{remark:r=0}
\normalfont In case the Mordell-Weil rank is $r=0$, we have $E(K)=E(K)_\text{tor}$. In particular, any two distinct sections are disjoint \cite[Cor. 8.30]{MWL}, hence every $k\geq 1$ is a gap number of $X$. For positive rank, the description of gap numbers is less trivial, thus our focus on $r\geq 1$.
\end{remark}

\subsection{Bounds $c_\text{max},c_\text{min}$ for the contribution term}\label{subsection:bounds}\
\indent We define the estimates $c_\text{max},c_\text{min}$ for the contribution term $\sum_v\text{contr}_v(P)$ and state some simple facts about them. We also provide an example to illustrate how they are computed.

\indent The need for these estimates comes from the following. Suppose we are given a section $P\in E(K)$ whose height $h(P)$ is known and we want to determine $P\cdot O$. In case $P\in E(K)^0$ we have a direct answer, namely $P\cdot O=h(P)/2-1$ by the height formula (\ref{height_formula_P}). However if $P\notin E(K)^0$, the computation of $P\cdot O$ depends on the contribution term $c_P:=\sum_{v\in R}\text{contr}_v(P)$, which by Table~\ref{contributions} depends on how $P$ intersects the reducible fibers of $\pi$. Usually we do not have this intersection data at hand, which is why we need estimates for $c_P$ not depending on $P$. 
\begin{defi} If the set $R$ of reducible fibers of $\pi$ is not empty, we define
\begin{align*}
c_\text{max}&:=\sum_{v\in R}\max\{\text{contr}_v(P)\mid P\in E(K)\},\\
c_\text{min}&:=\min\left\{\text{contr}_v(P)>0\mid P\in E(K), v\in R\right\}.
\end{align*}
\end{defi}

\begin{remark}
\normalfont The case $R=\emptyset$ only occurs when $X$ has Mordell-Weil rank $r=8$ (No. 1 in Table~\ref{tabela_completa}). In this case $E(K)^0=E(K)$ and $\sum_{v\in R}\text{contr}_v(P)=0$ $\forall P\in E(K)$, hence we adopt the convention $c_\text{max}=c_\text{min}=0$.
\end{remark}
\begin{remark}
\normalfont We use $c_\text{max},c_\text{min}$ as bounds for $c_P:=\sum_{v}\text{contr}_v(P)$. For our purposes it is not necessary to know whether $c_P$ actually attains one of these bounds for some $P$, so that $c_\text{max},c_\text{min}$ should be understood as hypothetical values.
\end{remark}

We state some facts about $c_\text{max},c_\text{min}$.

\begin{lemma}\label{lemma:bounds_are_actually_bounds}
Let $X$ be a rational elliptic surface with Mordell-Weil rank $r\geq 1$. If $\pi$ admits a reducible fiber, then:
\begin{enumerate}[i)]
\item $c_\text{min}>0$.
\item $c_\text{max}<4$.
\item ${\normalfont c_\text{min}\leq \sum_{v\in R}\text{contr}_v(P)\leq c_\text{max}}$ $\forall P\notin E(K)^0$. For $P\in E(K)^0$, only the second inequality holds.
\item If {\normalfont $\sum_{v\in R}\text{contr}_v(P)=c_\text{min}$}, then {\normalfont $\text{contr}_{v'}(P)=c_\text{min}$} for some $v'$ and {\normalfont $\text{contr}_v(P)=0$} for $v\neq v'$.
\end{enumerate}

\end{lemma}
\noindent\textit{Proof.} Item i) is immediate from the definition of $c_\text{min}$. For ii) it is enough to check the values of $c_\text{max}$ directly in Table~\ref{tabela_completa}. For iii), the second inequality follows from the definition of $c_\text{max}$ and clearly holds for any $P\in E(K)$. If $P\notin E(K)^0$, then $c_P:=\sum_v\text{contr}_v(P)>0$, so $\text{contr}_{v_0}(P)>0$ for some $v_0$. Therefore $c_P\geq \text{contr}_{v_0}(P)\geq c_\text{min}$.

For iv), let $\sum_v\text{contr}_v(P)=c_\text{min}$. Assume by contradiction that there are distinct $v_1,v_2$ such that $\text{contr}_{v_i}(P)>0$ for $i=1,2$. By definition of $c_\text{min}$ we have $c_\text{min}\leq\text{contr}_{v_i}(P)$ for $i=1,2$ so 
$$c_\text{min}=\sum_v\text{contr}_v(P)\geq \text{contr}_{v_1}(P)+\text{contr}_{v_2}(P)\geq 2c_\text{min},$$
\noindent which is absurd because $c_\text{min}>0$ by i). Therefore there is only one $v'$ with $\text{contr}_{v'}(P)>0$, while $\text{contr}_v(P)=0$ for all $v\neq v'$. In particular, $\text{contr}_{v'}(P)=c_\text{min}$.\fim
\\ \\
\noindent\textbf{Explicit computation.} Once we know the lattice $T$ associated with the reducible fibers of $\pi$ (Section~\ref{subsection:MW_lattice}), the computation of $c_\text{max},c_\text{min}$ is simple. For a fixed $v\in R$, the extreme values of the local contribution $\text{contr}_v(P)$ are given in Table~\ref{extreme_contributions}, which is derived from Table~\ref{contributions}. We provide an example to illustrate this computation. 

\begin{table}[h]
\begin{center}
\centering
\begin{tabular}{c|c|c} 
\hline
\multirow{2}{*}{\hfil $T_v$} & \multirow{2}{*}{\hfil $\max\{\text{contr}_v(P)\mid P\in E(K)\}$} & \multirow{2}{*}{\hfil $\min\{\text{contr}_v(P)>0\mid P\in E(K)\}$}\\ 
&\\
\hline
\multirow{2}{*}{\hfil $A_{n-1}$} & \multirow{2}{*}{\hfil $\frac{\ell(n-\ell)}{n}$, where $\ell:=\left\lfloor\frac{n}{2}\right\rfloor$} & \multirow{2}{*}{\hfil $\frac{n-1}{n}$}\\ 
&\\
\hline
\multirow{2}{*}{\hfil $D_{n+4}$} & \multirow{2}{*}{\hfil $1+\frac{n}{4}$} & \multirow{2}{*}{\hfil $1$}\\ 
&\\
\hline
\multirow{2}{*}{\hfil $E_6$} & \multirow{2}{*}{\hfil $\frac{4}{3}$} & \multirow{2}{*}{\hfil $\frac{4}{3}$}\\ 
&\\
\hline
\multirow{2}{*}{\hfil $E_7$} & \multirow{2}{*}{\hfil $\frac{3}{2}$} & \multirow{2}{*}{\hfil $\frac{3}{2}$}\\ 
&\\
\hline
\end{tabular}
\caption{Extreme values of $\text{contr}_v(P)$.}\label{extreme_contributions}
\end{center}
\end{table}

\noindent\textbf{Example:} Let $\pi$ with fiber configuration $(\text{I}_4,\text{IV},\text{III},\text{I}_1)$. The reducible fibers are $\text{I}_4,\text{IV},\text{III}$, so $T=A_3\oplus A_2\oplus A_1$. By Table~\ref{extreme_contributions}, the maximal contributions for $A_3,A_2,A_1$ are $\frac{2\cdot 2}{4}=1$, $\frac{2}{3}$, $\frac{1}{2}$ respectively. The minimal positive contributions are $\frac{1\cdot 3}{4}=\frac{3}{4}$, $\frac{2}{3}$, $\frac{1}{2}$ respectively. Then
\begin{align*}
c_\text{max}&=1+\frac{2}{3}+\frac{1}{2}=\frac{13}{6},\\
c_\text{min}&=\min\left\{\frac{3}{4},\frac{2}{3},\frac{1}{2}\right\}=\frac{1}{2}.
\end{align*}

\subsection{The difference $\Delta=c_\text{max}-c_\text{min}$}\label{subsection:presenting_Delta}\
\indent In this section we explain why the value of $\Delta:=c_\text{max}-c_\text{min}$ is relevant to our discussion, specially in Subsection~\ref{subsection:sufficient_conditions_Delta_leq_2}. We also verify that $\Delta<2$ in most cases and identify the exceptional ones in Table~\ref{table:Delta=2} and Table~\ref{table:Delta>2}.

\indent As noted in Subsection~\ref{subsection:bounds}, in case $P\notin E(K)^0$ and $h(P)$ is known, the difficulty of determining $P\cdot O$ lies in the contribution term $c_P:=\sum_{v\in R}\text{contr}_v(P)$. In particular, the range of possible values for $c_P$ determines the possibilities for $P\cdot O$. This range is measured by the difference
$$\Delta:=c_\text{max}-c_\text{min}.$$

Hence a smaller $\Delta$ means a better control over the intersection number $P\cdot O$, which is why $\Delta$ plays an important role in determining possible intersection numbers. In Subsection~\ref{subsection:necessary_sufficient_conditions_Delta_leq_2} we assume $\Delta\leq 2$ and state necessary and sufficient conditions for having a pair $P_1,P_2$ such that $P_1\cdot P_2=k$ for a given $k\geq 0$. If however $\Delta>2$, the existence of such a pair is not guaranteed a priori, so a case-by-case treatment is needed. Fortunately by Lemma~\ref{lemma:cases_where_Delta_geq_2} the case $\Delta>2$ is rare. 

\begin{lemma}\label{lemma:cases_where_Delta_geq_2}
Let $X$ be a rational elliptic surface with Mordell-Weil rank $r\geq 1$. The only cases with $\Delta=2$ and $\Delta>2$ are in Table \ref{table:Delta=2} and \ref{table:Delta>2} respectively. In particular we have $\Delta<2$ whenever $E(K)$ is torsion-free.

\begin{table}[h]
\begin{center}
\centering
\begin{tabular}{ccccc}
No. & $T$ & $E(K)$ & $c_\text{max}$ & $c_\text{min}$\\ 
\hline
\multirow{2}{*}{\hfil 24} & \multirow{2}{*}{\hfil $A_1^{\oplus 5}$} & \multirow{2}{*}{\hfil ${A_1^*}^{\oplus 3}\oplus\Bbb{Z}/2\Bbb{Z}$} & \multirow{2}{*}{\hfil $\frac{5}{2}$} & \multirow{2}{*}{\hfil $\frac{1}{2}$}\\
&\\
\multirow{2}{*}{\hfil 38} & \multirow{2}{*}{\hfil $A_3\oplus A_1^{\oplus 3}$} & \multirow{2}{*}{\hfil $A_1^*\oplus\langle 1/4\rangle\oplus\Bbb{Z}/2\Bbb{Z}$} & \multirow{2}{*}{\hfil $\frac{5}{2}$} & \multirow{2}{*}{\hfil $\frac{1}{2}$}\\
&\\
\multirow{2}{*}{\hfil 53} & \multirow{2}{*}{\hfil $A_5\oplus A_1^{\oplus 2}$} & \multirow{2}{*}{\hfil $\langle 1/6\rangle\oplus \Bbb{Z}/2\Bbb{Z}$} & \multirow{2}{*}{\hfil $\frac{5}{2}$} & \multirow{2}{*}{\hfil $\frac{1}{2}$}\\
&\\
\multirow{2}{*}{\hfil 57} & \multirow{2}{*}{\hfil $D_4\oplus A_1^{\oplus 3}$} & \multirow{2}{*}{\hfil $A_1^*\oplus(\Bbb{Z}/2\Bbb{Z})^{\oplus 2}$} & \multirow{2}{*}{\hfil $\frac{5}{2}$} & \multirow{2}{*}{\hfil $\frac{1}{2}$}\\
&\\
\multirow{2}{*}{\hfil 58} & \multirow{2}{*}{\hfil $A_3^{\oplus 2}\oplus A_1$} & \multirow{2}{*}{\hfil $A_1^*\oplus\Bbb{Z}/4\Bbb{Z}$} & \multirow{2}{*}{\hfil $\frac{5}{2}$} & \multirow{2}{*}{\hfil $\frac{1}{2}$} \\
&\\
\multirow{2}{*}{\hfil 61} & \multirow{2}{*}{\hfil $A_2^{\oplus 3}\oplus A_1$} & \multirow{2}{*}{\hfil  $\langle 1/6\rangle\oplus \Bbb{Z}/3\Bbb{Z}$} & \multirow{2}{*}{\hfil $\frac{5}{2}$} & \multirow{2}{*}{\hfil $\frac{1}{2}$}\\
&\\
\hline
\end{tabular}
\caption{Cases with $\Delta=2$}\label{table:Delta=2}
\end{center}
\end{table}
\newpage
\begin{table}[h]
\begin{center}
\centering
\begin{tabular}{cccccc}
No. & $T$ & $E(K)$ & $c_\text{max}$ & $c_\text{min}$ & $\Delta$\\ 
\hline
\multirow{3}{*}{\hfil 41} & \multirow{3}{*}{\hfil $A_2\oplus A_1^{\oplus 4}$} & \multirow{3}{*}{\hfil $\frac{1}{6}\left(\begin{matrix}2 & 1\\1 & 2\end{matrix}\right)\oplus\Bbb{Z}/2\Bbb{Z}$} & \multirow{3}{*}{\hfil $\frac{8}{3}$} & \multirow{3}{*}{\hfil $\frac{1}{2}$} & \multirow{3}{*}{\hfil $\frac{13}{6}$}\\
&\\
&\\
\multirow{2}{*}{\hfil 42} & \multirow{2}{*}{\hfil $A_1^{\oplus 6}$} & \multirow{2}{*}{\hfil ${A_1^*}^{\oplus 2}\oplus (\Bbb{Z}/2\Bbb{Z})^{\oplus 2}$} & \multirow{2}{*}{\hfil $3$} & \multirow{2}{*}{\hfil $\frac{1}{2}$} & \multirow{2}{*}{\hfil $\frac{5}{2}$}\\
&\\
\multirow{2}{*}{\hfil 59} & \multirow{2}{*}{\hfil $A_3\oplus A_2\oplus A_1^{\oplus 2}$} & \multirow{2}{*}{\hfil $\langle 1/12\rangle\oplus \Bbb{Z}/2\Bbb{Z}$} & \multirow{2}{*}{\hfil $\frac{8}{3}$} & \multirow{2}{*}{\hfil $\frac{1}{2}$} & \multirow{2}{*}{\hfil $\frac{13}{6}$}\\
&\\
\multirow{2}{*}{\hfil 60} & \multirow{2}{*}{\hfil $A_3\oplus A_1^{\oplus 4}$} & \multirow{2}{*}{\hfil  $\langle 1/4\rangle\oplus (\Bbb{Z}/2\Bbb{Z})^{\oplus 2}$} & \multirow{2}{*}{\hfil $3$} & \multirow{2}{*}{\hfil $\frac{1}{2}$} & \multirow{2}{*}{\hfil $\frac{5}{2}$}\\
&\\
\hline
\end{tabular}
\caption{Cases with $\Delta>2$}\label{table:Delta>2}
\end{center}
\end{table}
\end{lemma}
\noindent\textit{Proof.} By searching Table~\ref{tabela_completa} for all cases with $\Delta=2$ and $\Delta>2$, we obtain Table~\ref{table:Delta=2} and Table~\ref{table:Delta>2} respectively. Notice in particular that in both tables the torsion part of $E(K)$ is always nontrivial. Consequently, if $E(K)$ is torsion-free, then $\Delta<2$.\fim

\subsection{The quadratic form $Q_X$}\label{subsection:presenting_Q_X}\
\indent We define the positive-definite quadratic form with integer coefficients $Q_X$ derived from the height pairing. The relevance of $Q_X$ is due to the fact that some conditions for having $P_1\cdot P_2=k$ for some $P_1,P_2\in E(K)$ can be stated in terms of what integers can be represented by $Q_X$ (see Corollary~\ref{coro:necessary_conditions_Q_X} and Proposition~\ref{prop:summary_of_sufficient_conditions}).

The definition of $Q_X$ consists in clearing denominators of the rational quadratic form induced by the height pairing; the only question is how to find a scale factor that works in every case. More precisely, if $E(K)$ has rank $r\geq 1$ and $P_1,...,P_r$ are generators of its free part, then $q(x_1,...,x_r):=h(x_1P_1+...+x_rP_r)$ is a quadratic form with coefficients in $\Q$; we define $Q_X$ by multiplying $q$ by some integer $d>0$ so as to produce coefficients in $\Z$. We show that $d$ may always be chosen as the determinant of the narrow lattice $E(K)^0$.

\begin{defi}
Let $X$ with $r\geq 1$. Let $P_1,...,P_r$ be generators of the free part of $E(K)$. Define
$$Q_X(x_1,...,x_r):=(\det E(K)^0)\cdot h(x_1P_1+...+x_rP_r).$$
\end{defi}
We check that the matrix representing $Q_X$ has entries in $\Z$, therefore $Q_X$ has coefficients in $\Z$.

\begin{lemma}\label{lemma:adjugate}
Let $A$ be the matrix representing the quadratic form $Q_X$, i.e. $Q(x_1,...,x_r)=x^tAx$, where $x:=(x_1,...,x_r)^t$. Then $A$ has integer entries. In particular, $Q_X$ has integer coefficients.
% and is the adjugate of the Gram matrix $A^0$ of $E(K)^0$, i.e. $A^0\cdot A=(\det A^0)\cdot I_r$.
\end{lemma}
\noindent\textit{Proof.} Let $P_1,...,P_r$ be generators of the free part of $E(K)$ and let $L:=E(K)^0$. The free part of $E(K)$ is isomorphic to the dual lattice $L^*$ \cite[Main Thm.]{OguisoShioda}, so we may find generators $P_1^0,...,P_r^0$ of $L$ such that the Gram matrix $B^0:=(\langle P_i^0,P_j^0\rangle)_{i,j}$ of $L$ is the inverse of the Gram matrix $B:=(\langle P_i,P_j\rangle)_{i,j}$ of $L^*$.

\newpage
We claim that $Q_X$ is represented by the adjugate matrix of $B^0$, i.e. the matrix $\text{adj}(B^0)$ such that $B^0\cdot\text{adj}(B^0)=(\det B^0)\cdot I_r$, where $I_r$ is the $r\times r$ identity matrix. Indeed, by construction $B$ represents the quadratic form $h(x_1P_1+...+x_rP_r)$, therefore\begin{align*}
Q_X(x_1,...,x_r)&=(\det E(K)^0)\cdot h(x_1P_1+...+x_rP_r)\\
&=(\det B^0)\cdot x^tBx\\
&=(\det B^0)\cdot x^t(B^0)^{-1}x\\
&=x^t\text{adj}(B^0)x,
\end{align*}

as claimed. To prove that $A:=\text{adj}(B^0)$ has integer coefficients, notice that the Gram matrix $B^0$ of $L=E(K)^0$ has integer coefficients (as $E(K)^0$ is an even lattice), then so does $A$.\fim
\\ \\
\indent We close this subsection with a simple consequence of the definition of $Q_X$.

\begin{lemma}\label{lemma:Q_X_represents_dm} 
If $h(P)=m$ for some $P\in E(K)$, then $Q_X$ represents $d\cdot m$, where $d:=\det E(K)^0$.
\end{lemma}
\noindent\textit{Proof.} Let $P_1,...,P_r$ be generators for the free part of $E(K)$. Let $P=a_1P_1+...+a_rP_r+Q$, where $a_i\in\Z$ and $Q$ is a torsion element (possibly zero). Since torsion sections do not contribute to the height pairing, then $h(P-Q)=h(P)=m$. Hence 
\begin{align*}
Q_X(a_1,...,a_r)&=d\cdot h(a_1P_1+...+a_rP_r)\\
&=d\cdot h(P-Q)\\
&=d\cdot m.\text{\fim}
\end{align*}

\section{Intersection with a torsion section}\label{section:the_role_of_torsion_sections}\
\indent Before dealing with more technical details in Section~\ref{section:existence_of_a_pair}, we explain how torsion sections can be of help in our investigation, specially in Subsection~\ref{subsection:sufficient_conditions_Delta_leq_2}.

We first note some general properties of torsion sections. As the height pairing is positive-definite on $E(K)/E(K)_\text{tor}$, torsion sections are inert in the sense that for each $Q\in E(K)_\text{tor}$ we have $\langle Q,P\rangle=0$ for all $P\in E(K)$. Moreover, in the case of rational elliptic surfaces, torsion sections also happen to be mutually disjoint:

\begin{teor}{\normalfont \cite[Lemma 1.1]{MirandaPersson}}\label{thm:torsion_sections_are_disjoint}
On a rational elliptic surface, $Q_1\cdot Q_2=0$ for any distinct $Q_1,Q_2\in E(K)_\text{tor}$. In particular, if $O$ is the neutral section, then $Q\cdot O=0$ for all $Q\in E(K)_\text{tor}\setminus\{O\}$.
\end{teor}

\begin{remark}
\normalfont As stated in \cite[Lemma 1.1]{MirandaPersson}, Theorem~\ref{thm:torsion_sections_are_disjoint} holds for elliptic surfaces over $\Bbb{C}$ even without assuming $X$ is rational. However, for an arbitrary algebraically closed field the rationality hypothesis is needed, and a proof can be found in \cite[Cor. 8.30]{MWL}.
\end{remark}

By taking advantage of the properties above, we use torsion sections to help us find $P_1,P_2\in E(K)$ such that $P_1\cdot P_2=k$ for a given $k\in\Bbb{Z}_{\geq 0}$. This is particularly useful when $\Delta\geq 2$, in which case $E(K)_\text{tor}$ is not trivial by Lemma~\ref{lemma:cases_where_Delta_geq_2}. %As explained in Subsection~\ref{subsection:presenting_Delta}, the most difficult cases are the ones where $\Delta\geq 2$. These cases have the common property that $E(K)_\text{tor}$ is not trivial (Lemma~\ref{lemma:cases_where_Delta\geq_2}), hence the use of torison sections is particularly helpful.

The idea is as follows. Given $k\in\Bbb{Z}_{\geq 0}$, suppose we can find $P\in E(K)^0$ with height $h(P)=2k$. By the height formula (\ref{height_formula_P}), $P\cdot O=k-1<k$, which is not yet what we need. In the next lemma we show that replacing $O$ with a torsion section $Q\neq O$ gives $P\cdot Q=k$, as desired.
\newpage
\begin{lemma}\label{lemma:PO_plus_one}
Let $P\in E(K)^0$ such that $h(P)=2k$. Then $P\cdot Q=k$ for all $Q\in E(K)_\text{tor}\setminus\{O\}$.
\end{lemma}
\noindent\textit{Proof.} Assume there is some $Q\in E(K)_\text{tor}\setminus\{O\}$. By Theorem \ref{thm:torsion_sections_are_disjoint}, $Q\cdot O=0$ and by the height formula (\ref{height_formula_P}), $2k=2+2(P\cdot O)-0$, hence $P\cdot O=k-1$. We use the height formula (\ref{height_formula_PQ}) for $\langle P,Q\rangle$ in order to conclude that $P\cdot Q=k$. Since $P\in E(K)^0$, it intersects the neutral component $\Theta_{v,0}$ of every reducible fiber $\pi^{-1}(v)$, so $\text{contr}_v(P,Q)=0$ for all $v\in R$. Hence
\begin{align*}
0&=\langle P,Q\rangle\\
&=1+P\cdot O+Q\cdot O-P\cdot Q-\sum_{v\in R}\text{contr}_v(P,Q)\\
&=1+(k-1)+0-P\cdot Q-0\\
&=k-P\cdot Q.\text{\fim}
\end{align*}

\section{Existence of a pair of sections with a given intersection number}\label{section:existence_of_a_pair}\
\indent Given $k\in\Bbb{Z}_{\geq 0}$, we state necessary and (in most cases) sufficient conditions for having\linebreak $P_1\cdot P_2=k$ for some $P_1,P_2\in E(K)$. Necessary conditions are stated in generality in Subsection~\ref{subsection:necessary_conditions}, while sufficient ones depend on the value of $\Delta$ and are treated separately in Subsection~\ref{subsection:sufficient_conditions_Delta_leq_2}. In Subsection~\ref{subsection:summary_of_sufficient_conditions}, we collect all sufficient conditions proven in this section.

\subsection{Necessary Conditions}\label{subsection:necessary_conditions}\
\indent If $k\in\Bbb{Z}_{\geq 0}$, we state necessary conditions for having $P_1\cdot P_2=k$ for some sections $P_1,P_2\in E(K)$. We note that the value of $\Delta$ is not relevant in this subsection, although it plays a decisive role for sufficient conditions in Subsection~\ref{subsection:sufficient_conditions_Delta_leq_2}.

\begin{lemma}\label{lemma:necessary_conditions}
Let $k\in\Bbb{Z}_{\geq 0}$. If $P_1\cdot P_2=k$ for some $P_1,P_2\in E(K)$, then one of the following holds:
\begin{enumerate}[i)]
\item $h(P)=2+2k$ for some $P\in E(K)^0$.
\item $h(P)\in [2+2k-c_\text{max},\,2+2k-c_\text{min}]$ for some $P\notin E(K)^0$.
\end{enumerate}
\end{lemma}

\noindent\textit{Proof.} Without loss of generality we may assume $P_2$ is the neutral section, so that $P_1\cdot O=k$. By the height formula~(\ref{height_formula_P}), $h(P_1)=2+2k-c$, where $c:=\sum_v\text{contr}_v(P_1)$. If $P_1\in E(K)^0$, then $c=0$ and $h(P_1)=2+2k$, hence i) holds. If $P_1\notin E(K)^0$, then $c_\text{min}\leq c\leq c_\text{max}$ by Lemma~\ref{lemma:bounds_are_actually_bounds}. But $h(P_1)=2+2k-c$, therefore $2+2k-c_\text{max}\leq h(P_1)\leq 2+2k-c_\text{min}$, i.e. ii) holds.\fim

\begin{coro}\label{coro:necessary_conditions_Q_X}
Let $k\in\Bbb{Z}_{\geq 0}$. If $P_1\cdot P_2=k$ for some $P_1,P_2\in E(K)$, then $Q_X$ represents some integer in $[d\cdot (2+2k-c_\text{max}),d\cdot(2+2k)]$, where $d:=\det E(K)^0$.
\end{coro}
\noindent\textit{Proof.} We apply Lemma~\ref{lemma:necessary_conditions} and rephrase it in terms of $Q_X$. If i) holds, then $Q_X$ represents $d\cdot (2+2k)$ by Lemma~\ref{lemma:Q_X_represents_dm}. But if ii) holds, then $h(P)\in[2+2k-c_\text{max},2+2k-c_\text{min}]$ and by Lemma~\ref{lemma:Q_X_represents_dm}, $Q_X$ represents $d\cdot h(P)\in[d\cdot (2+2k-c_\text{max}),d\cdot(2+2k-c_\text{min})]$. In both i) and ii), $Q_X$ represents some integer in $[d\cdot(2+2k-c_\text{max}),d\cdot(2+2k)]$.\fim

\newpage

\subsection{Sufficient conditions when $\Delta\leq 2$}\label{subsection:sufficient_conditions_Delta_leq_2}\
\indent In this subsection we state sufficient conditions for having $P_1\cdot P_2=k$ for some $P_1,P_2\in E(K)$ under the assumption that $\Delta\leq 2$. By Lemma~\ref{lemma:cases_where_Delta_geq_2}, this covers almost all cases (more precisely, all but No. 41, 42, 59, 60 in Table~\ref{tabela_completa}). We treat $\Delta<2$ and $\Delta=2$ separately, as the latter needs more attention.

\subsubsection{The case $\Delta<2$}\label{subsubsection:the_case_Delta<2}\
\indent We first prove Lemma~\ref{lemma:sufficient_conditions_Delta<2}, which gives sufficient conditions assuming $\Delta<2$, then Corollary~\ref{coro:sufficient_conditions_Q_X_Delta<2}, which states sufficient conditions in terms of integers represented by $Q_X$. This is followed by Corollary~\ref{coro:sufficient_conditions_mu_Delta<2}, which is a simplified version of Corollary~\ref{coro:sufficient_conditions_Q_X_Delta<2}.

\begin{lemma}\label{lemma:sufficient_conditions_Delta<2}
Assume $\Delta<2$ and let $k\in\Bbb{Z}_{\geq 0}$. If $h(P)\in[2+2k-c_\text{max}, 2+2k-c_\text{min}]$ for some $P\notin E(K)^0$, then $P_1\cdot P_2=k$ for some $P_1,P_2\in E(K)$.
\end{lemma}
\noindent\textit{Proof.} Let $O\in E(K)$ be the neutral section. By the height formula (\ref{height_formula_P}), $h(P)=2+2(P\cdot O)-c$, where $c:=\sum_v\text{contr}_v(P)$. Since $h(P)\in [2+2k-c_\text{max},2+2k-c_\text{min}]$, then
\begin{align*}
2+2k-c_\text{max}&\leq 2+2(P\cdot O)-c\leq 2+2k-c_\text{min}\\
\Rightarrow \frac{c-c_\text{max}}{2} &\leq P\cdot O-k\leq \frac{c-c_\text{min}}{2}.
\end{align*}

Therefore $P\cdot O-k$ is an integer in $I:=\left[\frac{c-c_\text{max}}{2},\frac{c-c_\text{min}}{2}\right]$. We prove that $0$ is the only integer in $I$, so that $P\cdot O-k=0$, i.e. $P\cdot O=k$. First notice that $c\neq 0$, as $P\notin E(K)^0$. By Lemma~\ref{lemma:bounds_are_actually_bounds} iii), $c_\text{min}\leq c\leq c_\text{max}$, consequently $\frac{c-c_\text{max}}{2}\leq 0\leq \frac{c-c_\text{min}}{2}$, i.e. $0\in I$. Moreover $\Delta<2$ implies that $I$ has length $\frac{c_\text{max}-c_\text{min}}{2}=\frac{\Delta}{2}<1$, so $I$ contains no integer except $0$ as desired.\fim

\begin{remark}
\normalfont Lemma~\ref{lemma:sufficient_conditions_Delta<2} also applies when $c_\text{max}=c_\text{min}$, in which case the closed interval degenerates into a point.
\end{remark}

\indent The following corollary of Lemma~\ref{lemma:sufficient_conditions_Delta<2} states a sufficient condition in terms of integers represented by the quadratic form $Q_X$ (Section~\ref{subsection:presenting_Q_X}).

\begin{coro}\label{coro:sufficient_conditions_Q_X_Delta<2}
Assume $\Delta<2$ and let $d:=\det E(K)^0$. If $Q_X$ represents an integer not divisible by $d$ in the interval $[d\cdot(2+2k-c_\text{max}),d\cdot(2+2k-c_\text{min})]$, then $P_1\cdot P_2=k$ for some $P_1,P_2\in E(K)$.
\end{coro}
\textit{Proof.} Let $a_1,...,a_r\in\Bbb{Z}$ such that $Q_X(a_1,...,a_r)\in[d\cdot (2+2k-c_\text{max}),d\cdot(2+2k-c_\text{min})]$ with $d\nmid Q_X(a_1,...,a_r)$. Let $P:=a_1P_1+...+a_rP_r$, where $P_1,...,P_r$ are generators of the free part of $E(K)$. Then $d\nmid Q_X(a_1,...,a_r)=d\cdot h(P)$, which implies that $h(P)\notin\Bbb{Z}$. In particular $P\notin E(K)^0$ since $E(K)^0$ is an integer lattice. Moreover $h(P)=\frac{1}{d}Q_X(a_1,...,a_r)\in [2+2k-c_\text{max},2+2k-c_\text{min}]$ and we are done by Lemma~\ref{lemma:sufficient_conditions_Delta<2}.\fim

\newpage

\indent The next corollary, although weaker than Corollary~\ref{coro:sufficient_conditions_Q_X_Delta<2}, is more practical for concrete examples and is frequently used in Subsection~\ref{subsection:surfaces_with_a_1-gap}. It does not involve finding integers represented by $Q_X$, but only finding perfect squares in an interval depending on the minimal norm $\mu$ (Subsection~\ref{subsection:MW_lattice}).

\begin{coro}\label{coro:sufficient_conditions_mu_Delta<2}
Assume $\Delta<2$. If there is a perfect square $n^2\in\left[\frac{2+2k-c_\text{max}}{\mu},\frac{2+2k-c_\text{min}}{\mu}\right]$ such that $n^2\mu\notin\Bbb{Z}$, then $P_1\cdot P_2=k$ for some $P_1,P_2\in E(K)$.
\end{coro}
\noindent\textit{Proof.} Take $P\in E(K)$ such that $h(P)=\mu$. Since $h(nP)=n^2\mu\notin\Bbb{Z}$, we must have $nP\notin E(K)^0$ as $E(K)^0$ is an integer lattice. Moreover $h(nP)=n^2\mu\in[2+2k-c_\text{max},2+2k-c_\text{min}]$ and we are done by Lemma~\ref{lemma:sufficient_conditions_Delta<2}.\fim

\subsubsection{The case $\Delta=2$}\label{subsubsection:the_case_Delta=2}\
\indent The statement of sufficient conditions for $\Delta=2$ is almost identical to the one for $\Delta<2$: the only difference is that the closed interval Lemma~\ref{lemma:sufficient_conditions_Delta<2} is substituted by a right half-open interval in Lemma~\ref{lemma:sufficient_conditions_Delta=2}. This change, however, is associated with a technical difficulty in the case when a section has minimal contribution term, thus the separate treatment for $\Delta=2$.

The results are presented in the following order. First we prove Lemma~\ref{lemma:minimal_contribution_Delta=2}, which is a statement about sections whose contribution term is minimal. Next we prove Lemma~\ref{lemma:sufficient_conditions_Delta=2}, which states sufficient conditions for $\Delta=2$, then Corollaries~\ref{coro:sufficient_conditions_Q_X_Delta=2} and \ref{coro:sufficient_conditions_mu_Delta=2}.

\begin{lemma}\label{lemma:minimal_contribution_Delta=2}
Assume $\Delta=2$. If there is $P\in E(K)$ such that {\normalfont $\sum_{v\in R}\text{contr}_v(P)=c_\text{min}$}, then $P\cdot Q=P\cdot O+1$ for every $Q\in E(K)_\text{tor}\setminus\{O\}$.
\end{lemma}
\noindent\textit{Proof.} If $Q\in E(K)_\text{tor}\setminus\{O\}$, then $Q\cdot O=0$ by Theorem~\ref{thm:torsion_sections_are_disjoint}. Moreover, by the height formula (\ref{height_formula_PQ}), 
$$0=\langle P,Q\rangle=1+P\cdot O+0-P\cdot Q-\sum_{v\in R}\text{contr}_v(P,Q).\,\,(*)$$

Hence it suffices to show that $\text{contr}_v(P,Q)=0$ $\forall v\in R$. By Lemma~\ref{lemma:bounds_are_actually_bounds}~iv), $\text{contr}_{v'}(P)=c_\text{min}$ for some $v'$ and $\text{contr}_v(P)=0$ for all $v\neq v'$. In particular $P$ meets $\Theta_{v,0}$, hence $\text{contr}_v(P,Q)=0$ for all $v\neq v'$. Thus from $(*)$ we see that $\text{contr}_{v'}(P,Q)$ is an integer, which we prove is $0$.

We claim that $T_{v'}=A_1$, so that $\text{contr}_{v'}(P,Q)=0$ or $\frac{1}{2}$ by Table~\ref{contributions}. In this case, as $\text{contr}_{v'}(P,Q)$ is an integer, it must be $0$, and we are done. To see that $T_{v'}=A_1$ we analyse $\text{contr}_{v'}(P)$. Since $\Delta=2$, then $c_\text{min}=\frac{1}{2}$ by Table~\ref{table:Delta=2} and $\text{contr}_{v'}(P)=c_\text{min}=\frac{1}{2}$. By Table~\ref{contributions}, this only happens if $T_{v'}=A_{n-1}$ and $\frac{1}{2}=\frac{i(n-i)}{n}$ for some $0\leq i<n$. The only possibility is $i=1, n=2$ and $T_{v'}=A_1$.\fim
\ \\
%\\ \\
%\indent Before we prove the main results of this section, we still need another technical lemma, which states that if a section $P$ has height $h(P)=2+2k-c_\text{max}$, then $P\cdot Q=k$ for some section $Q$.
%\newpage
%\begin{lemma}\label{lemma_PQ}
%Assume $\Delta=2$. Let $k\geq 0$ be an integer and assume there is $P\in E(K)$ such that $h(P)=2+2k-c_\text{max}$. Then $P\cdot Q=k$ for some $Q\in E(K)_\text{tor}$.
%\end{lemma}
%\noindent\textit{Proof.} By the height formula (\ref{height_formula_P}), $h(P)=2+2(P\cdot O)-c$, where $c:=\sum_{v\in R}\text{contr}_v(P)$. Together with the hypothesis on $h(P)$, we obtain $c=2(P\cdot O-k)+c_\text{max}$ $(*)$. Since $\Delta=2$, then $c_\text{max}=\frac{5}{2}$ by Table~\ref{table:Delta=2}, hence from $(*)$ we conclude that $c\notin\Z$. In particular $c\neq 0$, so $c_\text{min}\leq c\leq c_\text{max}$ $(**)$ by Lemma~\ref{lemma:bounds_are_actually_bounds}. Combining $(*)$ and $(**)$ gives 
%$$0\leq 2(k-P\cdot O)\leq c_\text{max}-c_\text{min}=\Delta=2.$$

%Hence either $P\cdot O=k$ or $P\cdot O=k-1$. In the first case, we are done, so assume $P\cdot O=k-1$. By $(*)$ and the fact that $\Delta=2$, we have $c=-2+c_\text{max}=c_\text{min}$. This allows us to apply Lemma~\ref{lemma:minimal_contribution_Delta=2} and conclude that $P\cdot Q=(k-1)+1=k$ for every $Q\in E(K)_\text{tor}\setminus\{O\}$.\fim
%\\ \\
\indent With the aid of Lemma~\ref{lemma:minimal_contribution_Delta=2} we are able to state sufficient conditions for $\Delta=2$.
\newpage
\begin{lemma}\label{lemma:sufficient_conditions_Delta=2}
Assume $\Delta=2$ and let $k\in\Bbb{Z}_{\geq 0}$. If $h(P)\in[2+2k-c_\text{max}, 2+2k-c_\text{min})$ for some $P\notin E(K)^0$, then $P_1\cdot P_2=k$ for some $P_1,P_2\in E(K)$.
\end{lemma}
\noindent\textit{Proof.} Let $O\in E(K)$ be the neutral section. By the height formula (\ref{height_formula_P}), $h(P)=2+2(P\cdot O)-c$, where $c:=\sum_v\text{contr}_v(P)$. We repeat the arguments from Lemma~\ref{lemma:sufficient_conditions_Delta<2}, in this case with the right half-open interval, so that the hypothesis that $h(P)\in [2+2k-c_\text{max},2+2k-c_\text{min})$, implies that $P\cdot O-k$ is an integer in $I':=\left[\frac{c-c_\text{max}}{2},\frac{c-c_\text{min}}{2}\right)$. 

Since $I'$ is half-open with length $\frac{c_\text{max}-c_\text{min}}{2}=\frac{\Delta}{2}=1$, then $I'$ contains exactly one integer. If $0\in I'$, then $P\cdot O-k=0$, i.e. $P\cdot O=k$ and we are done. Hence we assume $0\notin I'$. 

We claim that $P\cdot O=k-1$. First, notice that if $c>c_\text{min}$, then the inequalities $c_\text{min}<c\leq c_\text{max}$ give $\frac{c-c_\text{max}}{2}\leq 0<\frac{c-c_\text{min}}{2}$, i.e. $0\in I'$, which is a contradiction. Hence $c=c_\text{min}$. Since $\Delta=2$, then $I'=[-1,0)$, whose only integer is $-1$. Thus $P\cdot O-k=-1$, i.e. $P\cdot O=k-1$, as claimed. 

Finally, let $Q\in E(K)_\text{tor}\setminus\{O\}$, so that $P\cdot Q=P\cdot O+1=k$ by Lemma~\ref{lemma:minimal_contribution_Delta=2} and we are done. We remark that $E(K)_\text{tor}$ is not trivial by Table~\ref{table:Delta=2}, therefore such $Q$ exists.\fim
\\ \\
\indent The following corollaries are analogues to Corollary~\ref{coro:sufficient_conditions_Q_X_Delta<2} and Corollary~\ref{coro:sufficient_conditions_mu_Delta<2} adapted to $\Delta=2$. Similarly to the case $\Delta<2$, Corollary~\ref{coro:sufficient_conditions_Q_X_Delta=2} is stronger than Corollary~\ref{coro:sufficient_conditions_mu_Delta=2}, although the latter is more practical for concrete examples. We remind the reader that $\mu$ denotes the minimal norm (Subsection~\ref{subsection:MW_lattice}).

\begin{coro}\label{coro:sufficient_conditions_Q_X_Delta=2}
Assume $\Delta=2$ and let $d:=\det E(K)^0$. If $Q_X$ represents an integer not divisible by $d$ in the interval $[d\cdot(2+2k-c_\text{max}),d\cdot(2+2k-c_\text{min}))$, then $P_1\cdot P_2=k$ for some $P_1,P_2\in E(K)$.
\end{coro}
\noindent\textit{Proof.} We repeat the arguments in Corollary~\ref{coro:sufficient_conditions_Q_X_Delta<2}, in this case with the half-open interval.\fim

\begin{coro}\label{coro:sufficient_conditions_mu_Delta=2}
Assume $\Delta=2$. If there is a perfect square $n^2\in\left[\frac{2+2k-c_\text{max}}{\mu},\frac{2+2k-c_\text{min}}{\mu}\right)$ such that $n^2\mu\notin\Bbb{Z}$, then $P_1\cdot P_2=k$ for some $P_1,P_2\in E(K)$.
\end{coro}
\noindent\textit{Proof.} We repeat the arguments in Corollary~\ref{coro:sufficient_conditions_mu_Delta<2}, in this case with the half-open interval.\fim

\subsection{Necessary and sufficient conditions for $\Delta\leq 2$}\label{subsection:necessary_sufficient_conditions_Delta_leq_2}\
\indent For completeness, we present a unified statement of necessary and sufficient conditions assuming $\Delta\leq 2$, which follows naturally from results in Subsections~\ref{subsection:necessary_conditions} and \ref{subsection:sufficient_conditions_Delta_leq_2}.

\begin{lemma}\label{lemma:necessary_sufficient_conditions_Delta_leq_2}
Assume $\Delta\leq 2$ and let $k\in\Bbb{Z}_{\geq 0}$. Then $P_1\cdot P_2=k$ for some $P_1,P_2\in E(K)$ if and only if one of the following holds:
\begin{enumerate}[i)]
\item $h(P)=2+2k$ for some $P\in E(K)^0$.
\item {\normalfont $h(P)\in[2+2k-c_\text{max},2+2k-c_\text{min})$} for some $P\notin E(K)^0$.
\item {\normalfont $h(P)=2+2k-c_\text{min}$} and {\normalfont $\sum_{v\in R}\text{contr}_v(P)=c_\text{min}$} for some $P\in E(K)$.
\end{enumerate}
\end{lemma}
\noindent\textit{Proof.} If i) or iii) holds, then $P\cdot O=k$ directly by the height formula (\ref{height_formula_P}). But if ii) holds, it suffices to to apply Lemma~\ref{lemma:sufficient_conditions_Delta<2} when $\Delta<2$ and by Lemma~\ref{lemma:sufficient_conditions_Delta=2} when $\Delta=2$.

Conversely, let $P_1\cdot P_2=k$. Without loss of generality, we may assume $P_2=O$, so that $P_1\cdot O=k$. By the height formula (\ref{height_formula_P}), $h(P_1)=2+2k-c$, where $c:=\sum_v\text{contr}_v(P_1)$. 

If $c=0$, then $P_1\in E(K)^0$ and $h(P_1)=2+2k$, so i) holds. Hence we let $c\neq 0$, i.e. $P_1\notin E(K)^0$, so that $c_\text{min}\leq c\leq c_\text{max}$ by Lemma~\ref{lemma:bounds_are_actually_bounds}. In case $c=c_\text{min}$, then $h(P_1)=2+2k-c_\text{min}$ and iii) holds. Otherwise $c_\text{min}<c\leq c_\text{max}$, which implies $2+2k-c_\text{max}\leq h(P_1)<2+2k-c_\text{min}$, so ii) holds.\fim

\subsection{Summary of sufficient conditions}\label{subsection:summary_of_sufficient_conditions}\
\indent For the sake of clarity, we summarize in a single proposition all sufficient conditions for having $P_1\cdot P_2=k$ for some $P_1,P_2\in E(K)$ proven in this section.

\begin{prop}\label{prop:summary_of_sufficient_conditions}
Let $k\in\Bbb{Z}_{\geq 0}$. If one of the following holds, then $P_1\cdot P_2=k$ for some $P_1,P_2\in E(K)$.
\begin{enumerate}[1)]
\item $h(P)=2+2k$ for some $P\in E(K)^0$.
\item $h(P)=2k$ for some $P\in E(K)^0$ and $E(K)_\text{tor}$ is not trivial.
\item $\Delta<2$ and there is a perfect square $n^2\in \left[\frac{2+2k-c_\text{max}}{\mu},\frac{2+2k-c_\text{min}}{\mu}\right]$ with $n^2\mu\notin\Bbb{Z}$, where $\mu$ is the minimal norm (Subsection~\ref{subsection:MW_lattice}). In case $\Delta=2$, consider the right half-open interval.
\item $\Delta<2$ and the quadratic form $Q_X$ represents an integer not divisible by $d:=\det E(K)^0$ in the interval $[d\cdot(2+2k-c_\text{max}),d\cdot(2+2k-c_\text{min})]$. In case $\Delta=2$, consider the right half-open interval.
\end{enumerate}
\end{prop}
\noindent\textit{Proof.} In 1) a height calculation gives $2+2k=h(P)=2+2(P\cdot O)-0$, so $P\cdot O=k$. For 2), we apply Lemma~\ref{lemma:PO_plus_one} to conclude that $P\cdot Q=k$ for any $Q\in E(K)_\text{tor}\setminus\{O\}$. In 3) we use Corollary~\ref{coro:sufficient_conditions_mu_Delta<2} when $\Delta<2$ and Corollary~\ref{coro:sufficient_conditions_mu_Delta=2} when $\Delta=2$. In 4), we apply Corollary~\ref{coro:sufficient_conditions_Q_X_Delta<2} if $\Delta<2$ and Corollary~\ref{coro:sufficient_conditions_Q_X_Delta=2} if $\Delta=2$.\fim

\section{Main Results}\label{section:main_results}\
\indent We prove the four main theorems of this paper, which are independent applications of the results from Section~\ref{section:existence_of_a_pair}. The first two are general attempts to describe when and how gap numbers occur: Theorem~\ref{thm:gap_free_r_geq_5} tells us that large Mordell-Weil groups prevent the existence of gaps numbers, more precisely for Mordell-Weil rank $r\geq 5$; in Theorem~\ref{thm:gaps_probability_1_r=1,2} we show that for small Mordell-Weil rank, more precisely when $r\leq 2$, then gap numbers occur with probability $1$. The last two theorems, on the other hand, deal with explicit values of gap numbers: Theorem~\ref{thm:identification_of_gaps_r=1} provides a complete description of gap numbers in certain cases, while Theorem~\ref{thm:surfaces_with_a_1-gap} is a classification of cases with a $1$-gap.

\subsection{No gap numbers in rank $r\geq 5$}\label{subsection:gap_free_r_geq_5}\
\indent We show that if $E(K)$ has rank $r\geq 5$, then $X$ is gap-free. Our strategy is to prove that for every $k\in\Bbb{Z}_{\geq 0}$ there is some $P\in E(K)^0$ such that $h(P)=2+2k$, and by Proposition~\ref{prop:summary_of_sufficient_conditions} 1) we are done. We accomplish this in two steps. First we show that this holds when there is an embedding of $A_1^{\oplus}$ or of $A_4$ in $E(K)^0$ (Lemma~\ref{lemma:sublattice_A1_times_4_or_A4}). Second, we show that if $r\geq 5$, then such embedding exists, hence $X$ is gap-free (Theorem~\ref{thm:gap_free_r_geq_5}).
\newpage
\begin{lemma}\label{lemma:sublattice_A1_times_4_or_A4}
Assume $E(K)^0$ has a sublattice isomorphic to $A_1^{\oplus 4}$ or $A_4$. Then for every $\ell\in\Bbb{Z}_{\geq 0}$ there is $P\in E(K)^0$ such that $h(P)=2\ell$.
\end{lemma}
\noindent\textit{Proof.} First assume $A_1^{\oplus 4}\subset E(K)^0$ and let $P_1,P_2,P_3,P_4$ be generators for each factor $A_1$ in $A_1^{\oplus 4}$. Then $h(P_i)=2$ and $\langle P_i,P_j\rangle=0$ for distinct $i,j=1,2,3,4$. By Lagrange's four-square theorem \cite[\S 20.5]{HardyWright} there are integers $a_1,a_2,a_3,a_4$ such that $a_1^2+a_2^2+a_3^2+a_4^2=\ell$. Defining $P:=a_1P_1+a_2P_2+a_3P_3+a_4P_4\in A_1^{\oplus 4}\subset E(K)^0$, we have
$$h(P)=2a_1^2+2a_2^2+2a_3^2+2a_4^2=2\ell.$$

Now let $A_4\subset E(K)^0$ with generators $P_1,P_2,P_3,P_4$. Then $h(P_i)=2$ for $i=1,2,3,4$ and $\langle P_i,P_{i+1}\rangle=-1$ for $i=1,2,3$. We need to find integers $x_1,...,x_4$ such that $h(P)=2\ell$, where $P:=x_1P_1+...+x_4P_4\in A_4\subset E(K)^0$. Equivalently, we need that
$$\ell=\frac{1}{2}\langle P,P\rangle=x_1^2+x_2^2+x_3^2+x_4^2-x_1x_2-x_2x_3-x_3x_4.$$

Therefore $\ell$ must be represented by $q(x_1,...,x_4):=x_1^2+x_2^2+x_3^2+x_4^2-x_1x_2-x_2x_3-x_3x_4$. We prove that $q$ represents all positive integers. Notice that $q$ is positive-definite, since it is induced by $\langle\cdot,\cdot\rangle$. By Bhargava-Hanke's 290-theorem \cite{BhargavaHanke}[Thm. 1], $q$ represents all positive integers if and only if it represents the following integers:
$$2, 3, 5, 6, 7, 10, 13, 14, 15, 17, 19, 21, 22, 23, 26,$$
$$29, 30, 31, 34, 35, 37, 42, 58, 93, 110, 145, 203, 290.$$

The representation for each of the above is found in Table~\ref{table:representation_critical_integers}.\fim
\\ \\
\indent We now prove the main theorem of this section.

\begin{teor}\label{thm:gap_free_r_geq_5} 
If $r\geq 5$, then $X$ is gap-free. 
\end{teor}
\noindent\textit{Proof.} We show that for every $k\geq 0$ there is $P\in E(K)^0$ such that $h(P)=2+2k$, so that by Proposition~\ref{prop:summary_of_sufficient_conditions} 1) we are done. Using Lemma~\ref{lemma:sublattice_A1_times_4_or_A4} it suffices to prove that $E(K)^0$ has a sublattice isomorphic to $A_1^{\oplus 4}$ or $A_4$.

The cases with $r\geq 5$ are No. 1-7 (Table~\ref{tabela_completa}). In No. 1-6, $E(K)^0= E_8,E_7,E_6,D_6,D_5,A_5$ respectively. Each of these admit an $A_4$ sublattice \cite[Lemmas 4.2,4.3]{Nishiyama}. In No. 7 we claim that $E(K)^0=D_4\oplus A_1$ has an $A_1^{\oplus 4}$ sublattice. This is the case because $D_4$ admits an $A_1^{\oplus 4}$ sublattice \cite[Lemma 4.5 (iii)]{Nishiyama}.\fim

%Indeed, let $P_1,P_2,P_3,P_4$ be generators for $D_4$ and $P_5$ a generator for $A_1$. Then $L:=\Bbb{Z}\langle P_1,P_3,P_4,P_5\rangle$ is clearly isomorphic to $A_1^{\oplus 4}$ (Figure \ref{D_4_A_1}).\fim
%\begin{figure}[h!]
%\begin{center}
%\includegraphics[scale=0.7]{D_4_A_1}
%\caption{Generators $P_1,P_3,P_4,P_5$ for the sublattice $L\simeq A_1^{\oplus 4}\subset D_4\oplus A_1$.}\label{D_4_A_1}
%\end{center}
%\end{figure}

\newpage
\begin{table}[h]
\begin{center}
\centering
\begin{tabular}{c|c} 
$n$ & $x_1,x_2,x_3,x_4$ with $x_1^2+x_2^2+x_3^2+x_4^2-x_1x_2-x_2x_3-x_3x_4=n$\\
\hline
$1$ & $1,0,0,0$\\
$2$ & $1,0,1,0$\\
$3$ & $1,1,2,0$\\
$5$ & $1,0,2,0$\\
$6$ & $1,1,-2,-1$\\
$7$ & $1,1,-2,0$\\
$10$ & $1,0,3,0$\\
$13$ & $2,0,3,0$\\
$14$ & $1,2,5,1$\\
$15$ & $1,5,5,2$\\
$17$ & $1,0,4,0$\\
$19$ & $1,5,3,-1$\\
$21$ & $1,5,0,0$\\
$22$ & $1,5,0,-1$\\
$23$ & $1,6,6,2$\\
$26$ & $1,0,5,0$\\
$29$ & $2,0,5,0$\\
$30$ & $1,5,0,-3$\\
$31$ & $1,3,-4,-2$\\
$34$ & $3,0,5,0$\\
$35$ & $1,2,-2,4$\\
$37$ & $1,0,6,0$\\
$42$ & $1,1,-4,3$\\
$58$ & $3,0,7,0$\\
$93$ & $1,1,-10,0$\\
$110$ & $1,-2,3,-8$\\
$145$ & $1,0,12,0$\\
$203$ & $1,-5,-9,8$\\
$290$ & $1,0,17,0$\\
\end{tabular}\caption{Representation of the critical integers in Bhargava-Hanke's 290-theorem.}\label{table:representation_critical_integers}
\end{center}
\end{table}

\subsection{Gaps with probability $1$ in rank $r\leq 2$}\label{subsection:gaps_probability_1_r=1,2}\
\indent Fix a rational elliptic surface $\pi:X\to\P^1$ with Mordell-Weil rank $r\leq 2$. We prove that if $k$ is a uniformly random natural number, then $k$ is a gap number with probability $1$. More precisely, if $G:=\{k\in\Bbb{N}\mid k\text{ is a gap number of }X\}$ is the set of gap numbers, then $G\subset\Bbb{N}$ has density $1$, i.e. 
$$d(G):=\lim_{n\to\infty}\frac{\#G\cap\{1,...,n\}}{n}=1.$$

\newpage

\indent We adopt the following strategy. If $k\in\Bbb{N}\setminus G$, then $P_1\cdot P_2=k$ for some $P_1,P_2\in E(K)$ and by Corollary~\ref{coro:necessary_conditions_Q_X} the quadratic form $Q_X$ represents some integer $t$ depending on $k$. This defines a function $\Bbb{N}\setminus G\to T$, where $T$ is the set of integers represented by $Q_X$. Since $Q_X$ is a quadratic form on $r\leq 2$ variables, $T$ has density $0$ in $\Bbb{N}$ by Lemma~\ref{lemma:representable_integers_density_0}. By analyzing the pre-images of $\Bbb{N}\setminus G\to T$, in Theorem~\ref{thm:gaps_probability_1_r=1,2} we conclude that $d(\Bbb{N}\setminus G)=d(T)=0$, hence $d(G)=1$ as desired.

\begin{lemma}\label{lemma:representable_integers_density_0}
Let $Q$ be a positive-definite quadratic form on $r=1,2$ variables with integer coefficients. Then the set of integers represented by $Q$ has density $0$ in $\Bbb{N}$.
\end{lemma}
\noindent\textit{Proof.} Let $S$ be the set of integers represented by $Q$. If $d$ is the greatest common divisor of the coefficients of $Q$, let $S'$ be the set of integers representable by the primitive form $Q':=\frac{1}{d}\cdot Q$. By construction $S'$ is a rescaling of $S$, so $d(S)=0$ if and only if $d(S')=0$.

If $r=1$, then $Q'(x_1)=x_1^2$ and $S'$ is the set of perfect squares, so clearly $d(S')=0$. If $r=2$, then $Q'$ is a binary quadratic form and the number  of elements in $S'$ bounded from above by $x>0$ is given by $C\cdot \frac{x}{\sqrt{\log x}}+o(x)$ with $C>0$ a constant and $\lim_{x\to\infty}\frac{o(x)}{x}=0$ \cite[p. 91]{Bernays}. Thus 
$$d(S')=\lim_{x\to\infty}\frac{C}{\sqrt{\log x}}+\frac{o(x)}{x}=0.\text{\fim}$$

We now prove the main result of this section.
\begin{teor}\label{thm:gaps_probability_1_r=1,2}
Let $\pi:X\to\P^1$ be a rational elliptic surface with Mordell-Weil rank $r\leq 2$. Then the set $G:=\{k\in\Bbb{N}\mid k\text{ is a gap number of }X\}$ of gap numbers of $X$ has density $1$ in $\Bbb{N}$.
\end{teor}
\noindent\textit{Proof.} If $r=0$, then the claim is trivial by Remark~\ref{remark:r=0}, hence we may assume $r=1,2$. We prove that $S:=\Bbb{N}\setminus G$ has density $0$. If $S$ is finite, there is nothing to prove. Otherwise, let $k_1<k_2<...$ be the increasing sequence of all elements of $S$. By Corollary~\ref{coro:necessary_conditions_Q_X}, for each $n$ there is some $t_n\in J_{k_n}:=[d\cdot (2+2k_n-c_\text{max}),d\cdot(2+2k_n)]$ represented by the quadratic form $Q_X$. Let $T$ be the set of integers represented by $Q_X$ and define the function $f\colon \Bbb{N}\setminus G\to T$ by $k_n\mapsto t_n$. Since $Q_X$ has $r=1,2$ variables, $T$ has density $0$ by Lemma~\ref{lemma:representable_integers_density_0}.

For $N>0$, let $S_N:=S\cap\{1,...,N\}$ and $T_N:=T\cap\{1,...,N\}$. Since $T$ has density zero, $\#T_N=o(N)$, i.e. $\frac{\# T_N}{N}\to 0$ when $N\to\infty$ and we need to prove that $\#S_N=o(N)$. We analyze the function $f$ restricted to $S_N$. Notice that as $t_n\in J_{k_n}$, then $k_n\leq N$ implies $t_n\leq d\cdot (2+2k_n)\leq d\cdot(2+2N)$. Hence the restriction $g:=f|_{S_N}$ can be regarded as a function $g:S_N\to T_{d\cdot(2+2k)}$. 

We claim that $\#g^{-1}(t)\leq 2$ for all $t\in T_{d\cdot(2+2N)}$, in which case $\#S_N\leq 2\cdot\#T_{d\cdot(2+2N)}=o(N)$ and we are done. Assume by contradiction that $g^{-1}(t)$ contains three distinct elements, say $k_{\ell_1}<k_{\ell_2}<k_{\ell_3}$ with $t=t_{\ell_1}=t_{\ell_2}=t_{\ell_3}$. Since $t_{\ell_i}\in J_{k_{\ell_i}}$ for each $i=1,2,3$, then $t\in J_{k_{\ell_1}}\cap J_{k_{\ell_2}}\cap J_{k_{\ell_3}}$. We prove that $J_{k_{\ell_1}}$ and $J_{k_{\ell_3}}$ are disjoint, which yields a contradiction. Indeed, since $k_{\ell_1}<k_{\ell_2}<k_{\ell_3}$, in particular $k_{\ell_3}-k_{\ell_1}\geq 2$, therefore $d\cdot(2+2k_{\ell_1})\leq d\cdot(2+2k_{\ell_3}-4)$. But $c_\text{max}<4$ by Lemma~\ref{lemma:bounds_are_actually_bounds}, so $d\cdot(2+2k_{\ell_1})<d\cdot(2+2k_{\ell_3}-c_\text{max})$, i.e. $\max J_{k_{\ell_1}}<\min J_{k_{\ell_3}}$. Thus $J_{k_{\ell_1}}\cap J_{k_{\ell_3}}=\emptyset$, as desired.\fim

\subsection{Identification of gaps when $E(K)$ is torsion-free with rank $r=1$}\label{subsection:identification_of_gaps_r=1}\
\indent The results in Subsections~\ref{subsection:gap_free_r_geq_5} and \ref{subsection:gaps_probability_1_r=1,2} concern the existence and the distribution of gap numbers. In the following subsections we turn our attention to finding gap numbers explicitly. In this subsection we give a complete description of gap numbers assuming $E(K)$ is torsion-free with rank $r=1$. Such descriptions are difficult in the general case, but our assumption guarantees that each $E(K),E(K)^0$ is generated by a single element and that $\Delta<2$ by Lemma \ref{lemma:cases_where_Delta_geq_2}, which makes the problem more accessible. 

%This section is organized as follows. In Lemma~\ref{lemma:generators_r=1} we collect some straight-forward facts about generators of $E(K),E(K)^0$ when $E(K)$ has rank $r=1$. As an application we prove Corollary~\ref{E0=A1}, which is particularly useful in Section~\ref{RES_gap_2}. In Lemma~\ref{lemma:necessary_sufficient_conditions_r=1} we assume $E(K)$ is torsion-free with rank $r=1$ and state, for a given $k\geq 0$, necessary and sufficient conditions for having $P\cdot P'=k$ for some $P,P'\in E(K)$. At last we prove Theorem~\ref{thm:identification_of_gaps_r=1}, which is the main result of this section.

We organize this subsection as follows. First we point out some trivial facts about generators of $E(K),E(K)^0$ when $r=1$ in Lemma~\ref{lemma:generators_r=1}. Next we state necessary and sufficient conditions for having $P_1\cdot P_2=k$ when $E(K)$ is torsion-free with $r=1$ in Lemma~\ref{lemma:necessary_sufficient_conditions_r=1}. As an application of the latter, we prove Theorem~\ref{thm:identification_of_gaps_r=1}, which is the main result of the subsection.

\begin{lemma}\label{lemma:generators_r=1}
Let $X$ be a rational elliptic surface with Mordell-Weil rank $r=1$. If $P$ generates the free part of $E(K)$, then 
\begin{enumerate}[a)]
\item $h(P)=\mu$.
\item $1/\mu$ is an even integer.
\item $E(K)^0$ is generated by $P_0:=(1/\mu)P$ and $h(P_0)=1/\mu$.
\end{enumerate}
\end{lemma}
\noindent\textit{Proof.} Item  a) is clear. Items b), c) follow from the fact that $E(K)^0$ is an even lattice and that $E(K)\simeq L^*\oplus E(K)_\text{tor}$, where $L:=E(K)^0$ \cite[Main Thm.]{OguisoShioda}.\fim 
\ \\
\indent In what follows we use Lemma~\ref{lemma:generators_r=1} and results from Section~\ref{section:existence_of_a_pair} to state necessary and sufficient conditions for having $P_1\cdot P_2=k$ for some $P_1,P_2\in E(K)$ in case $E(K)$ is torsion-free with $r=1$. 

\begin{lemma}\label{lemma:necessary_sufficient_conditions_r=1}
Assume $E(K)$ is torsion-free with rank $r=1$. Then $P_1\cdot P_2=k$ for some $P_1,P_2\in E(K)$ if and only if one of the following holds:
\begin{enumerate}[i)]
\item $\mu\cdot(2+2k)$ is a perfect square.
\item There is a perfect square $n^2\in\left[\frac{2+2k-c_\text{max}}{\mu},\frac{2+2k-c_\text{min}}{\mu}\right]$ such that $\mu\cdot n\notin\Bbb{Z}$.

%By Table~\ref{tabela_completa}, $E(K)$ is torsion-free with $r=1$ precisely in \normalfont{No. = 43, 45-47, 49, 50, 55, 56}.
\end{enumerate}
\end{lemma}
\noindent\textit{Proof.} By Lemma~\ref{lemma:generators_r=1}, $E(K)$ is generated by some $P$ with $h(P)=\mu$ and $E(K)^0$ is generated by $P_0:=n_0P$, where $n_0:=\frac{1}{\mu}\in 2\Bbb{Z}$.

First assume that $P_1\cdot P_2=k$ for some $P_1,P_2$. Without loss of generality we may assume $P_2=O$. Let $P_1=nP$ for some $n\in\Bbb{Z}$. We show that $P_1\in E(K)^0$ implies i) while $P_1\notin E(K)^0$ implies ii). 

If $P_1\in E(K)^0$, then $n_0\mid n$, hence $P_1=nP=mP_0$, where $m:=\frac{n}{n_0}$. By the height formula (\ref{height_formula_P}), $2+2k=h(P_1)=h(mP_0)=m^2\cdot \frac{1}{\mu}$. Hence $\mu\cdot (2+2k)=m^2$, i.e. i) holds. 

If $P_1\notin E(K)^0$, then $n_0\nmid n$, hence $\mu\cdot n=\frac{n}{n_0}\notin\Bbb{Z}$. Moreover, $h(P_1)=n^2h(P)=n^2\mu$ and by the height formula (\ref{height_formula_P}), $n^2\mu=h(P)=2+2k-c$, where $c:=\sum_v\text{contr}_v(P_1)\neq 0$. The inequalities $c_\text{min}\leq c\leq c_\text{max}$ then give $\frac{2+2k-c_\text{max}}{\mu}\leq n^2\leq \frac{2+2k-c_\text{min}}{\mu}$. Hence ii) holds.

Conversely, assume i) or ii) holds. Since $E(K)$ is torsion-free, $\Delta<2$ by Lemma~\ref{lemma:cases_where_Delta_geq_2}, so we may apply Lemma~\ref{lemma:sufficient_conditions_Delta<2}. If i) holds, then $\mu\cdot(2+2k)=m^2$ for some $m\in\Bbb{Z}$. Since $mP_0\in E(K)^0$ and $h(mP_0)=\frac{m^2}{\mu}=2+2k$, we are done by Lemma \ref{lemma:sufficient_conditions_Delta<2} i). If ii) holds, the condition $\mu\cdot n\notin\Bbb{Z}$ is equivalent to $n_0\nmid n$, hence $nP\notin E(K)^0$. Moreover $n^2\in\left[\frac{2+2k-c_\text{max}}{\mu},\frac{2+2k-c_\text{min}}{\mu}\right]$, implies $h(nP)=n^2\mu\in[2+2k-c_\text{max},2+2k-c_\text{min}]$. By Lemma~\ref{lemma:sufficient_conditions_Delta<2} ii), we are done.\fim
\\ \\
\indent By applying Lemma~\ref{lemma:necessary_sufficient_conditions_r=1} to all possible cases where $E(K)$ is torsion-free with rank $r=1$, we obtain the main result of this subsection.
\newpage
\begin{teor}\label{thm:identification_of_gaps_r=1}
If $E(K)$ is torsion-free with rank $r=1$, then all the gap numbers of $X$ are described in Table~\ref{table:description_gaps_r=1}.
\begin{table}[h]
\begin{center}
\centering
\begin{tabular}{cccc} 
\hline
\multirow{2}{*}{No.} & \multirow{2}{*}{$T$} & \multirow{2}{*}{$\begin{matrix}k\text{ is a gap number}\Leftrightarrow \text{none of}\\ \text{the following are perfect squares}\end{matrix}$} & \multirow{2}{*}{first gap numbers}\\ 
& \\
\hline
\multirow{2}{*}{43} & \multirow{2}{*}{$E_7$} & \multirow{2}{*}{$k+1$, $4k+1$} & \multirow{2}{*}{$1,4$}\\ 
& \\
\hline
\multirow{2}{*}{45} & \multirow{2}{*}{$A_7$} & \multirow{2}{*}{$\frac{k+1}{4}$, $16k,...,16k+9$} & \multirow{2}{*}{$8,11$}\\ 
& \\
\hline
\multirow{2}{*}{46} & \multirow{2}{*}{$D_7$} & \multirow{2}{*}{$\frac{k+1}{2}$, $8k+1,...,8k+4$} & \multirow{2}{*}{$2,5$}\\ 
& \\
\hline
\multirow{2}{*}{47} & \multirow{2}{*}{$A_6\oplus A_1$} & \multirow{2}{*}{$\frac{k+1}{7}$, $28k-3,...,28k+21$} & \multirow{2}{*}{$12,16$}\\ 
& \\
\hline
\multirow{2}{*}{49} & \multirow{2}{*}{$E_6\oplus A_1$} & \multirow{2}{*}{$\frac{k+1}{3}$, $12k+1,...,12k+9$} & \multirow{2}{*}{$3,7$}\\ 
& \\
\hline
\multirow{2}{*}{50} & \multirow{2}{*}{$D_5\oplus A_2$} & \multirow{2}{*}{$\frac{k+1}{6}$, $24k+1,...,24k+16$} & \multirow{2}{*}{$6,11$}\\ 
& \\
\hline
\multirow{2}{*}{55} & \multirow{2}{*}{$A_4\oplus A_3$} & \multirow{2}{*}{$\frac{k+1}{10}$, $40k-4,...,40k+25$} & \multirow{2}{*}{$16,20$}\\ 
& \\
\hline
\multirow{2}{*}{56} & \multirow{2}{*}{$A_4\oplus A_2\oplus A_1$} & \multirow{2}{*}{$\frac{k+1}{15}$, $60k-11,...,60k+45$} & \multirow{2}{*}{$22,27$}\\ 
& \\
\hline
\end{tabular}\caption{Description of gap numbers when $E(K)$ is torsion-free with $r=1$.}\label{table:description_gaps_r=1}
\end{center}
\end{table}
\end{teor}
\noindent\textit{Proof.} For the sake of brevity we restrict ourselves to No. 55. The other cases are treated similarly. Here $c_\text{max}=\frac{2\cdot 3}{5}+\frac{2\cdot 2}{4}=\frac{11}{5}$, $c_\text{min}=\min\left\{\frac{4}{5},\frac{3}{4}\right\}=\frac{3}{4}$ and $\mu=1/20$. 

By Lemma~\ref{lemma:necessary_sufficient_conditions_r=1}, $k$ is a gap number if and only if neither i) nor ii) occurs. Condition i) is that $\frac{2+2k}{20}=\frac{k+1}{10}$ is a perfect square. Condition ii) is that $\left[\frac{2+2k-c_\text{max}}{\mu},\frac{2+2k-c_\text{min}}{\mu}\right]=[40k-4,40k+25]$ contains some $n^2$ with $20\nmid n$. We check that $20\nmid n$ for every $n$ such that $n^2=40k+\ell$, with $\ell=-4,...,25$. Indeed, if $20\mid n$, then $400\mid n^2$ and in particular $40\mid n^2$. Then $40\mid (n^2-40k)=\ell$, which is absurd.\fim

\subsection{Surfaces with a $1$-gap}\label{subsection:surfaces_with_a_1-gap}\
\indent In Subsection~\ref{subsection:identification_of_gaps_r=1} we take each case in Table~\ref{table:description_gaps_r=1} and describe all its gap numbers. In this subsection we do the opposite, which is to fix a number and describe all cases having it as a gap number. We remind the reader that our motivating problem (Section~\ref{section:intro}) was to determine when there are sections $P_1,P_2$ such that $P_1\cdot P_2=1$, which induce a conic bundle having $P_1+P_2$ as a reducible fiber. The answer for this question is the main theorem of this subsection:

\begin{teor}\label{thm:surfaces_with_a_1-gap}
Let $X$ be a rational elliptic surface. Then $X$ has a $1$-gap if and only if $r=0$ or $r=1$ and $\pi$ has a {\normalfont $\text{III}^*$} fiber.
\end{teor}
\newpage
Our strategy for the proof is the following. We already know that a $1$-gap exists whenever $r=0$ (Theorem~\ref{thm:torsion_sections_are_disjoint}) or when $r=1$ and $\pi$ has a $\text{III}^*$ fiber (Theorem~\ref{thm:identification_of_gaps_r=1}, No. 43). Conversely, we need to find $P_1,P_2$ with $P_1\cdot P_2=1$ in all cases with $r\geq 1$ and $T\neq E_7$. 

First we introduce two lemmas, which solve most cases with little computation, and leave the remaining ones for the proof of Theorem~\ref{thm:surfaces_with_a_1-gap}. In both Lemma~\ref{lemma:E0_height_4} and Lemma~\ref{lemma:E0=An} our goal is to analyze the narrow lattice $E(K)^0$ and apply Proposition~\ref{prop:summary_of_sufficient_conditions} to detect cases without a $1$-gap.

\begin{lemma}\label{lemma:E0_height_4}
If one of the following holds, then $h(P)=4$ for some $P\in E(K)^0$.
\begin{enumerate}[a)]
\item The Gram matrix of $E(K)^0$ has a $4$ in its main diagonal.
\item There is an embedding of $A_n\oplus A_m$ in $E(K)^0$ for some $n,m\geq 1$.
\item There is an embedding of $A_n,D_n$ or $E_n$ in $E(K)^0$ for some $n\geq 3$.
\end{enumerate}

\end{lemma}
\noindent\textit{Proof.} Case a) is trivial. Assuming b), we take generators $P_1,P_2$ from $A_n,A_m$ respectively with $h(P_1)=h(P_2)=2$. Since $A_n,A_m$ are in direct sum, $\langle P_1,P_2\rangle=0$, hence $h(P_1+P_2)=4$, as desired. If c) holds, then the fact that $n\geq 3$ allows us to choose two elements $P_1,P_2$ among the generators of $L_1=A_n,D_n$ or $E_n$ such that $h(P_1)=h(P_2)=2$ and $\langle P_1,P_2\rangle=0$. Thus $h(P_1+P_2)=4$ as claimed.\fim
\begin{coro}\label{a_b_c_coro}
In the following cases, $X$ does not have a $1$-gap.
\begin{itemize}
\item $r\geq 3:$ all cases except possibly {\normalfont No. 20}.
\item $r=1,2:$ cases {\normalfont No. 25, 26, 30, 32-36, 38, 41, 42, 46, 52, 54, 60}.
\end{itemize}
\end{coro}
\noindent\textit{Proof.} We look at column $E(K)^0$ in Table~\ref{tabela_completa} to find which cases satisfy one of the conditions a), b), c) from Lemma \ref{lemma:E0_height_4}. 
\begin{enumerate}[a)]
\item Applies to No. 12, 17, 19, 22, 23, 25, 30, 32, 33, 36, 38, 41, 46, 52, 54, 60.
\item Applies to No. 10, 11, 14, 15, 18, 24, 26, 34, 35, 42.
\item Applies to No. 1-10, 13, 16, 21.
\end{enumerate}

In particular, this covers all cases with $r\geq 3$ (No. 1-24) except No. 20. By Lemma~\ref{lemma:E0_height_4} in each of these cases there is $P\in E(K)^0$ with $h(P)=4$ and we are done by Proposition~\ref{prop:summary_of_sufficient_conditions} 1).\fim
\\ \\
\indent In the next lemma we also analyze $E(K)^0$ to detect surfaces without a $1$-gap.

\begin{lemma}\label{lemma:E0=An}
Assume $E(K)^0\simeq A_n$ for some $n\geq 1$ and that $E(K)$ has nontrivial torsion part. Then $X$ does not have a $1$-gap. This applies to cases {\normalfont No. 28, 39, 44, 48, 51, 57, 58} in Table~\ref{tabela_completa}.
\end{lemma}
\noindent\textit{Proof.} Take a generator $P$ of $E(K)^0$ with $h(P)=2$ and apply Proposition~\ref{prop:summary_of_sufficient_conditions} 2).\fim
\newpage
\indent We are ready to prove the main result of this subsection.
\\ \\
\noindent\textit{Proof of Theorem \ref{thm:surfaces_with_a_1-gap}.} We need to show that in all cases where $r\geq 1$ and $T\neq E_7$ there are $P_1,P_2\in E(K)$ such that $P_1\cdot P_2=1$. This corresponds to cases No. 1-61 except 43 in Table~\ref{tabela_completa}. 

The cases where $r=1$ and $E(K)$ is torsion-free can be solved by Theorem~\ref{a_b_c_coro}, namely No. 45-47, 49, 50, 55, 56. Adding these cases to the ones treated in Corollary~\ref{a_b_c_coro} and Lemma~\ref{lemma:E0=An}, we have therefore solved the following:
$$\text{No. }1\text{-}19,\,21\text{-}26,\,28,\,30,\,32\text{-}36,\,38,\,39,\,41\text{-}52,\,54\text{-}58,\,60.$$

%We follow the method described in Section \ref{subsection:summary_of_sufficient_conditions}.

%1) We look for elements $P\in E(K)^0$ such that $h(P)=4$. This is done in Lemma~\ref{lemma:E0_height_4} and Corollary~\ref{a_b_c_coro} and solves the following cases:
%$$\text{No. = }1\text{-}19,\,21\text{-}26,\,30,\,32\text{-}36,\,38,\,41,\,42,\,46,\,52,\,54,\,60.$$

%2) We look for elements $P\in E(K)^0$ such that $h(P)=2$ when $E(K)$ has nontrivial torsion part. This is done in Lemma~\ref{lemma:E0=An} and applies to the following cases:
%$$\text{No. = }28,\,39,\,44,\,48,\,51,\,57,\,58.$$

%3) We use Theorem~\ref{thm:identification_of_gaps_r=1} with $k=1$. This shows that $X$ has a $1$-gap when $\pi$ has a $\text{III}^*$ fiber (No. $43$) and that $X$ does not have a $1$-gap in the following cases:
%$$\text{No.}=45\text{-}47,\,49,\,50,\,55,\,56.$$

For the remaining cases, we apply Proposition~\ref{prop:summary_of_sufficient_conditions} 3), which involves finding perfect squares in the interval $\left[\frac{4-c_\text{max}}{\mu},\frac{4-c_\text{min}}{\mu}\right]$ (see Table~\ref{table:perfect_squares_in_the_interval}), considering the half-open interval in the cases with $\Delta=2$ (No. 53, 61).

\begin{table}[h]
\begin{center}
\centering
\begin{tabular}{cccccc} 

\multirow{2}{*}{No.} & \multirow{2}{*}{$T$} & \multirow{2}{*}{$E(K)$} & \multirow{2}{*}{$\mu$} & \multirow{2}{*}{$I$} & \multirow{2}{*}{$n^2\in I$}\\ 
& \\
\hline
\multirow{2}{*}{\hfil 20} & \multirow{2}{*}{\hfil $A_2^{\oplus 2}\oplus A_1$} & \multirow{2}{*}{\hfil $A_2^*\oplus\langle1/6\rangle$} & \multirow{2}{*}{$\frac{1}{6}$} & \multirow{2}{*}{$[13,23]$} & \multirow{2}{*}{$4^2$}\\ 
&\\
\multirow{2}{*}{\hfil 27} & \multirow{2}{*}{\hfil $E_6$} & \multirow{2}{*}{\hfil $A_2^*$} & \multirow{2}{*}{$\frac{2}{3}$} & \multirow{2}{*}{$[4,4]$} & \multirow{2}{*}{$2^2$}\\ 
&\\
\multirow{2}{*}{\hfil 29} & \multirow{2}{*}{\hfil $A_5\oplus A_1$} & \multirow{2}{*}{\hfil $A_1^*\oplus\langle 1/6\rangle$} & \multirow{2}{*}{$\frac{1}{6}$} & \multirow{2}{*}{$[12,21]$} & \multirow{2}{*}{$4^2$}\\
& \\ 
\multirow{3}{*}{\hfil 31} & \multirow{3}{*}{\hfil $A_4\oplus A_2$} & \multirow{3}{*}{\hfil $\frac{1}{15}\left(\begin{matrix} 2 & 1\\1 & 8\end{matrix}\right)$} & \multirow{3}{*}{$\frac{2}{15}$} & \multirow{3}{*}{$[16,21]$} & \multirow{3}{*}{$4^2$}\\
& \\
& \\
\multirow{2}{*}{\hfil 37} & \multirow{2}{*}{\hfil $A_3\oplus A_2\oplus A_1$} & \multirow{2}{*}{\hfil $A_1^*\oplus\langle 1/12\rangle$} & \multirow{2}{*}{$\frac{1}{12}$} & \multirow{2}{*}{$[22,28]$} & \multirow{2}{*}{$5^2$}\\
& \\ 
\multirow{2}{*}{\hfil 40} & \multirow{2}{*}{\hfil $A_2^{\oplus 2}\oplus A_1^{\oplus 2}$} & \multirow{2}{*}{\hfil $\langle 1/6\rangle^{\oplus 2}$}  & \multirow{2}{*}{$\frac{1}{6}$} & \multirow{2}{*}{$\left[10,21\right]$} & \multirow{2}{*}{$4^2$}\\
& \\ 
\multirow{2}{*}{\hfil 53} & \multirow{2}{*}{\hfil $A_5\oplus A_1^{\oplus 2}$} & \multirow{2}{*}{\hfil $\langle 1/6\rangle\oplus \Bbb{Z}/2\Bbb{Z}$} &  \multirow{2}{*}{$\frac{1}{6}$} & \multirow{2}{*}{$[9,12]$} & \multirow{2}{*}{$3^2$}\\
& \\ 
\multirow{2}{*}{\hfil 59} & \multirow{2}{*}{\hfil $A_3\oplus A_2\oplus A_1^{\oplus 2}$} & \multirow{2}{*}{\hfil $\langle 1/12\rangle\oplus \Bbb{Z}/2\Bbb{Z}$} & \multirow{2}{*}{$\frac{1}{12}$} & \multirow{2}{*}{$[16,42]$} & \multirow{2}{*}{$4^2,5^2,6^2$}\\
& \\ 
\multirow{3}{*}{\hfil 61} & \multirow{3}{*}{\hfil $A_2^{\oplus 3}\oplus A_1$} & \multirow{3}{*}{\hfil $\langle 1/6\rangle\oplus\Bbb{Z}/3\Bbb{Z}$} & \multirow{3}{*}{$\frac{1}{6}$} & \multirow{3}{*}{$[9,12]$} & \multirow{3}{*}{$3^2$}\\
& \\ 
\hline
\end{tabular}\caption{Perfect squares in the interval $I:=\left[\frac{4-c_\text{max}}{\mu},\frac{4-c_\text{min}}{\mu}\right]$.}\label{table:perfect_squares_in_the_interval}
\end{center}
\end{table}

\indent In No. 59 we have $\Delta>2$, so a particular treatment is needed. Let $T=T_{v_1}\oplus T_{v_2}\oplus T_{v_3}\oplus T_{v_4}=A_3\oplus A_2\oplus A_1\oplus A_1$. If $P$ generates the free part of $E(K)$ and $Q$ generates its torsion part, then $h(P)=\frac{1}{12}$ and $4P+Q$ meets the reducible fibers at $\Theta_{v_1,2},\Theta_{v_2,1},\Theta_{v_3,1},\Theta_{v_4,1}$ \cite{Kurumadani}[Example 1.7]. By Table~\ref{contributions} and the height formula (\ref{height_formula_P}), 
$$\frac{4^2}{12}=h(4P+Q)=2+2(4P+Q)\cdot O-\frac{2\cdot 2}{4}-\frac{1\cdot 2}{3}-\frac{1}{2}-\frac{1}{2},$$

hence $(4P+Q)\cdot O=1$, as desired.\fim

\newpage

\section{Appendix}\label{section:appendix}\
\indent We reproduce part of the table in \cite[Main Th.]{OguisoShioda} with data on Mordell-Weil lattices of rational elliptic surfaces with Mordell-Weil rank $r\geq 1$. We only add columns $c_\text{max},c_\text{min},\Delta$.

\begin{center}
\begin{longtable}{cccccccc}\label{tabela_completa}
\multirow{1}{*}{No.} & \multirow{1}{*}{$r$} & \multirow{1}{*}{$T$} & \multirow{1}{*}{$E(K)^0$} & \multirow{1}{*}{$E(K)$} & \multirow{1}{*}{$c_\text{max}$} & \multirow{1}{*}{$c_\text{min}$} & \multirow{1}{*}{$\Delta$}\\
\hline
\multirow{1}{*}{1} & \multirow{1}{*}{$8$} & \multirow{1}{*}{$0$} & \multirow{1}{*}{$E_8$} & \multirow{1}{*}{$E_8$} & \multirow{1}{*}{0} & \multirow{1}{*}{0} & \multirow{1}{*}{0}\\
\hline
\multirow{2}{*}{2} & \multirow{2}{*}{$7$} & \multirow{2}{*}{$A_1$} & \multirow{2}{*}{$E_7$} & \multirow{2}{*}{$E_8^*$} & \multirow{2}{*}{$\frac{1}{2}$} & \multirow{2}{*}{$\frac{1}{2}$} & \multirow{2}{*}{$0$}\\
&\\
\hline
\multirow{2}{*}{3} & \multirow{2}{*}{$6$} & \multirow{2}{*}{$A_2$} & \multirow{2}{*}{$E_6$} & \multirow{2}{*}{$E_6^*$} & \multirow{2}{*}{$\frac{2}{3}$} & \multirow{2}{*}{$\frac{2}{3}$} & \multirow{2}{*}{$0$}\\
&\\
\multirow{2}{*}{4} & \multirow{2}{*}{} & \multirow{2}{*}{$A_1^{\oplus 2}$} & \multirow{2}{*}{$D_6$} & \multirow{2}{*}{$D_6^*$} & \multirow{2}{*}{$\frac{3}{2}$} & \multirow{2}{*}{$1$} & \multirow{2}{*}{$\frac{1}{2}$}\\
&\\
\hline
\multirow{2}{*}{5} & \multirow{2}{*}{$5$} & \multirow{2}{*}{$A_3$} & \multirow{2}{*}{$D_5$} & \multirow{2}{*}{$D_5^*$} & \multirow{2}{*}{$1$} & \multirow{2}{*}{$\frac{3}{4}$} & \multirow{2}{*}{$\frac{1}{4}$}\\
&\\
\multirow{2}{*}{6} & \multirow{2}{*}{} & \multirow{2}{*}{$A_2\oplus A_1$} & \multirow{2}{*}{$A_5$} & \multirow{2}{*}{$A_5^*$} & \multirow{2}{*}{$\frac{7}{6}$} & \multirow{2}{*}{$\frac{1}{2}$} & \multirow{2}{*}{$\frac{2}{3}$}\\
&\\
\multirow{2}{*}{7} & \multirow{2}{*}{} & \multirow{2}{*}{$A_1^{\oplus 3}$} & \multirow{2}{*}{$D_4\oplus A_1$} & \multirow{2}{*}{$D_4^*\oplus A_1^*$} & \multirow{2}{*}{$\frac{3}{2}$} & \multirow{2}{*}{$\frac{1}{2}$} & \multirow{2}{*}{$1$}\\
&\\
\hline
\multirow{2}{*}{8} & \multirow{2}{*}{$4$} & \multirow{2}{*}{$A_4$} & \multirow{2}{*}{$A_4$} & \multirow{2}{*}{$A_4^*$} & \multirow{2}{*}{$\frac{6}{5}$} & \multirow{2}{*}{$\frac{4}{5}$} & \multirow{2}{*}{$\frac{2}{5}$}\\
&\\
\multirow{2}{*}{9} & \multirow{2}{*}{} & \multirow{2}{*}{$D_4$} & \multirow{2}{*}{$D_4$} & \multirow{2}{*}{$D_4^*$} & \multirow{2}{*}{$1$} & \multirow{2}{*}{$1$} & \multirow{2}{*}{$0$}\\
&\\
\multirow{2}{*}{10} & \multirow{2}{*}{} & \multirow{2}{*}{$A_3\oplus A_1$} & \multirow{2}{*}{$A_3\oplus A_1$} & \multirow{2}{*}{$A_3^*\oplus A_1^*$} & \multirow{2}{*}{$\frac{3}{2}$} & \multirow{2}{*}{$\frac{1}{2}$} & \multirow{2}{*}{$1$}\\
&\\
\multirow{2}{*}{\hfil 11} & \multirow{2}{*}{\hfil } & \multirow{2}{*}{\hfil $A_2^{\oplus 2}$} & \multirow{2}{*}{$A_2^{\oplus 2}$} & \multirow{2}{*}{${A_2^*}^{\oplus 2}$} & \multirow{2}{*}{\hfil $\frac{4}{3}$} & \multirow{2}{*}{\hfil $\frac{2}{3}$} & \multirow{2}{*}{\hfil $\frac{2}{3}$}\\
&\\
\multirow{4}{*}{\hfil 12} & \multirow{4}{*}{} & \multirow{4}{*}{\hfil $A_2\oplus A_1^{\oplus 2}$} & \multirow{4}{*}{$\left(\begin{matrix}4 & -1 & 0 & 1\\-1 & 2 & -1 & 0\\0 & -1 & 2 & -1\\1 & 0 & -1 & 2\end{matrix}\right)$} & \multirow{4}{*}{$\frac{1}{6}\left(\begin{matrix}2 & 1 & 0 & -1\\1 & 5 & 3 & 1\\0 & 3 & 6 & 3\\-1 & 1 & 3 & 5\end{matrix}\right)$} & \multirow{4}{*}{\hfil $\frac{5}{3}$} & \multirow{4}{*}{\hfil $\frac{1}{2}$} & \multirow{4}{*}{\hfil $\frac{7}{6}$}\\
&\\
&\\
&\\
\multirow{2}{*}{\hfil 13} & \multirow{2}{*}{} & \multirow{2}{*}{\hfil $A_1^{\oplus 4}$} & \multirow{2}{*}{$D_4$} & \multirow{2}{*}{$D_4^*\oplus\Bbb{Z}/2\Bbb{Z}$} & \multirow{2}{*}{\hfil $2$} & \multirow{2}{*}{\hfil $\frac{1}{2}$} & \multirow{2}{*}{\hfil $\frac{3}{2}$}\\
&\\
\multirow{2}{*}{\hfil 14} & \multirow{2}{*}{} & \multirow{2}{*}{\hfil $A_1^{\oplus 4}$} & \multirow{2}{*}{$A_1^{\oplus 4}$} & \multirow{2}{*}{${A_1^*}^{\oplus 4}$} & \multirow{2}{*}{\hfil $2$} & \multirow{2}{*}{\hfil $\frac{1}{2}$} & \multirow{2}{*}{\hfil $\frac{3}{2}$}\\
&\\
\hline
\multirow{2}{*}{\hfil 15} & \multirow{2}{*}{$3$} & \multirow{2}{*}{\hfil $A_5$} & \multirow{2}{*}{$A_2\oplus A_1$} & \multirow{2}{*}{$A_2^*\oplus A_1^*$} & \multirow{2}{*}{\hfil $\frac{3}{2}$} & \multirow{2}{*}{\hfil $\frac{5}{6}$} & \multirow{2}{*}{\hfil $\frac{2}{3}$}\\
&\\
\multirow{2}{*}{\hfil 16} & \multirow{2}{*}{} & \multirow{2}{*}{$D_5$} & \multirow{2}{*}{$A_3$} & \multirow{2}{*}{$A_3^*$} & \multirow{2}{*}{\hfil $\frac{5}{4}$} & \multirow{2}{*}{\hfil $1$} & \multirow{2}{*}{\hfil $\frac{1}{4}$}\\
&\\
\multirow{3}{*}{17} & \multirow{3}{*}{} & \multirow{3}{*}{$A_4\oplus A_1$} & \multirow{3}{*}{$\left(\begin{matrix}4 & -1 & 1\\-1 & 2 & -1\\1 & -1 & 2\end{matrix}\right)$} & \multirow{3}{*}{$\frac{1}{10}\left(\begin{matrix}3 & 1 & -1\\1 & 7 & 3\\-1 & 3 & 7\end{matrix}\right)$} & \multirow{3}{*}{\hfil $\frac{17}{10}$} & \multirow{3}{*}{\hfil $\frac{1}{2}$} & \multirow{3}{*}{\hfil $\frac{6}{5}$}\\
&\\
&\\
\multirow{2}{*}{18} & \multirow{2}{*}{} & \multirow{2}{*}{$D_4\oplus A_1$} & \multirow{2}{*}{$A_1^{\oplus 3}$} & \multirow{2}{*}{${A_1^*}^{\oplus 3}$} & \multirow{2}{*}{\hfil $\frac{3}{2}$} & \multirow{2}{*}{\hfil $\frac{1}{2}$} & \multirow{2}{*}{\hfil $1$}\\
&\\
\multirow{3}{*}{19} & \multirow{3}{*}{} & \multirow{3}{*}{$A_3\oplus A_2$} & \multirow{3}{*}{$\left(\begin{matrix}2 & 0 & -1\\0 & 2 & -1\\-1 & -1 & 4\end{matrix}\right)$} & \multirow{3}{*}{$\frac{1}{12}\left(\begin{matrix}7 & 1 & 2\\1 & 7 & 2\\2 & 2 & 4\end{matrix}\right)$} & \multirow{3}{*}{\hfil $\frac{5}{3}$} & \multirow{3}{*}{\hfil $\frac{2}{3}$} & \multirow{3}{*}{\hfil $1$}\\
&\\
&\\
\multirow{2}{*}{20} & \multirow{2}{*}{} & \multirow{2}{*}{$A_2^{\oplus 2}\oplus A_1$} & \multirow{2}{*}{$A_2\oplus\langle 6\rangle$} & \multirow{2}{*}{$A_2^*\oplus\langle 1/6\rangle$} & \multirow{2}{*}{\hfil $\frac{11}{6}$} & \multirow{2}{*}{\hfil $\frac{1}{2}$} & \multirow{2}{*}{\hfil $\frac{4}{3}$}\\
&\\
\multirow{2}{*}{21} & \multirow{2}{*}{} & \multirow{2}{*}{$A_3\oplus A_1^{\oplus 2}$} & \multirow{2}{*}{$A_3$} & \multirow{2}{*}{$A_3^*\oplus\Z/2\Z$} & \multirow{2}{*}{$2$} & \multirow{2}{*}{$\frac{1}{2}$} & \multirow{2}{*}{$\frac{3}{2}$}\\
&\\
\multirow{2}{*}{22} & \multirow{2}{*}{} & \multirow{2}{*}{$A_3\oplus A_1^{\oplus 2}$} & \multirow{2}{*}{$A_1\oplus\langle 4\rangle$} & \multirow{2}{*}{$A_1^*\oplus\langle 1/4\rangle$} & \multirow{2}{*}{$2$} & \multirow{2}{*}{$\frac{1}{2}$} & \multirow{2}{*}{$\frac{3}{2}$}\\
&\\
\multirow{2}{*}{23} & \multirow{2}{*}{} & \multirow{2}{*}{$A_2\oplus A_1^{\oplus 3}$} & \multirow{2}{*}{$A_1\oplus\left(\begin{matrix}4 & -2\\-2 & 4\end{matrix}\right)$} & \multirow{2}{*}{$A_1^*\oplus\frac{1}{6}\left(\begin{matrix}2 & 1\\1 & 2\end{matrix}\right)$} & \multirow{2}{*}{$\frac{13}{6}$} & \multirow{2}{*}{$\frac{1}{2}$} & \multirow{2}{*}{$\frac{5}{3}$}\\
&\\
\multirow{2}{*}{24} & \multirow{2}{*}{} & \multirow{2}{*}{$A_1^{\oplus 5}$} & \multirow{2}{*}{$A_1^{\oplus 3}$} & \multirow{2}{*}{${A_1^*}^{\oplus 3}\oplus\Z/2\Z$} & \multirow{2}{*}{$\frac{5}{2}$} & \multirow{2}{*}{$\frac{1}{2}$} & \multirow{2}{*}{$2$}\\
&\\
\hline
\multirow{2}{*}{25} & \multirow{2}{*}{$2$} & \multirow{2}{*}{$A_6$} & \multirow{2}{*}{$\left(\begin{matrix}4 & -1\\-1& 2\end{matrix}\right)$} & \multirow{2}{*}{$\frac{1}{7}\left(\begin{matrix}2 & 1\\1 & 4\end{matrix}\right)$} & \multirow{2}{*}{$\frac{12}{7}$} & \multirow{2}{*}{$\frac{6}{7}$} & \multirow{2}{*}{$\frac{6}{7}$}\\
&\\
\multirow{2}{*}{26} & \multirow{2}{*}{} & \multirow{2}{*}{$D_6$} & \multirow{2}{*}{$A_1^{\oplus 2}$} & \multirow{2}{*}{${A_1^*}^{\oplus 2}$} & \multirow{2}{*}{$\frac{3}{2}$} & \multirow{2}{*}{$1$} & \multirow{2}{*}{$\frac{1}{2}$}\\
&\\
\multirow{2}{*}{27} & \multirow{2}{*}{} & \multirow{2}{*}{$E_6$} & \multirow{2}{*}{$A_2$} & \multirow{2}{*}{$A_2^*$} & \multirow{2}{*}{$\frac{4}{3}$} & \multirow{2}{*}{$\frac{4}{3}$} & \multirow{2}{*}{$0$}\\
&\\
\multirow{2}{*}{28} & \multirow{2}{*}{} & \multirow{2}{*}{$A_5\oplus A_1$} & \multirow{2}{*}{$A_2$} & \multirow{2}{*}{$A_2^*\oplus\Z/2\Z$} & \multirow{2}{*}{$2$} & \multirow{2}{*}{$\frac{1}{2}$} & \multirow{2}{*}{$\frac{3}{2}$}\\
&\\
\multirow{2}{*}{29} & \multirow{2}{*}{} & \multirow{2}{*}{$A_5\oplus A_1$} & \multirow{2}{*}{$A_1\oplus\langle 6\rangle$} & \multirow{2}{*}{$A_1^*\oplus\langle 1/6\rangle$} & \multirow{2}{*}{$2$} & \multirow{2}{*}{$\frac{1}{2}$} & \multirow{2}{*}{$\frac{3}{2}$}\\
&\\
\multirow{2}{*}{30} & \multirow{2}{*}{} & \multirow{2}{*}{$D_5\oplus A_1$} & \multirow{2}{*}{$A_1\oplus\langle 4\rangle$} & \multirow{2}{*}{$A_1^*\oplus\langle 1/4\rangle$} & \multirow{2}{*}{$\frac{7}{4}$} & \multirow{2}{*}{$\frac{1}{2}$} & \multirow{2}{*}{$\frac{5}{4}$}\\
&\\
\multirow{3}{*}{31} & \multirow{3}{*}{} & \multirow{3}{*}{$A_4\oplus A_2$} & \multirow{3}{*}{$\left(\begin{matrix}8 & -1\\-1 & 2\end{matrix}\right)$} & \multirow{3}{*}{$\frac{1}{15}\left(\begin{matrix}2 & 1\\1 & 8\end{matrix}\right)$} & \multirow{3}{*}{$\frac{28}{15}$} & \multirow{3}{*}{$\frac{2}{3}$} & \multirow{3}{*}{$\frac{6}{5}$}\\
&\\
&\\
\multirow{3}{*}{32} & \multirow{3}{*}{} & \multirow{3}{*}{$D_4\oplus A_2$} & \multirow{3}{*}{$\left(\begin{matrix}4 & -2\\-2 & 4\end{matrix}\right)$} & \multirow{3}{*}{$\frac{1}{6}\left(\begin{matrix}2 & 1\\1 & 2\end{matrix}\right)$} & \multirow{3}{*}{$\frac{5}{3}$} & \multirow{3}{*}{$\frac{2}{3}$} & \multirow{3}{*}{$1$}\\
&\\
&\\
\multirow{2}{*}{33} & \multirow{2}{*}{} & \multirow{2}{*}{$A_4\oplus A_1^{\oplus 2}$} & \multirow{2}{*}{$\left(\begin{matrix}6 & -2\\-2 & 4\end{matrix}\right)$} & \multirow{2}{*}{$\frac{1}{10}\left(\begin{matrix}2 & 1\\1 & 3\end{matrix}\right)$} & \multirow{2}{*}{$\frac{11}{5}$} & \multirow{2}{*}{$\frac{1}{2}$} & \multirow{2}{*}{$\frac{17}{10}$}\\
&\\
\multirow{2}{*}{34} & \multirow{2}{*}{} & \multirow{2}{*}{$D_4\oplus A_1^{\oplus 2}$} & \multirow{2}{*}{$A_1^{\oplus 2}$} & \multirow{2}{*}{${A_1^*}^{\oplus 2}$} & \multirow{2}{*}{$2$} & \multirow{2}{*}{$\frac{1}{2}$} & \multirow{2}{*}{$\frac{3}{2}$}\\
&\\
\multirow{2}{*}{35} & \multirow{2}{*}{} & \multirow{2}{*}{$A_3^{\oplus 2}$} & \multirow{2}{*}{$A_1^{\oplus 2}$} & \multirow{2}{*}{${A_1^*}^{\oplus 2}\oplus \Z/2\Z$} & \multirow{2}{*}{$2$} & \multirow{2}{*}{$\frac{3}{4}$} & \multirow{2}{*}{$\frac{5}{4}$}\\
&\\
\multirow{2}{*}{36} & \multirow{2}{*}{} & \multirow{2}{*}{$A_3^{\oplus 2}$} & \multirow{2}{*}{$\langle 4\rangle^{\oplus 2}$} & \multirow{2}{*}{$\langle 1/4\rangle^{\oplus 2}$} & \multirow{2}{*}{$2$} & \multirow{2}{*}{$\frac{3}{4}$} & \multirow{2}{*}{$\frac{5}{4}$}\\
&\\
\multirow{2}{*}{37} & \multirow{2}{*}{} & \multirow{2}{*}{$A_3\oplus A_2\oplus A_1$} & \multirow{2}{*}{$A_1\oplus\langle 12\rangle$} & \multirow{2}{*}{$A_1^*\oplus\langle 1/12\rangle$} & \multirow{2}{*}{$\frac{13}{6}$} & \multirow{2}{*}{$\frac{1}{2}$} & \multirow{2}{*}{$\frac{5}{3}$}\\
&\\
\multirow{2}{*}{38} & \multirow{2}{*}{} & \multirow{2}{*}{$A_3\oplus A_1^{\oplus 3}$} & \multirow{2}{*}{$A_1\oplus\langle 4\rangle$} & \multirow{2}{*}{$A_1^*\oplus\langle 1/4\rangle\oplus\Z/2\Z$} & \multirow{2}{*}{$\frac{5}{2}$} & \multirow{2}{*}{$\frac{1}{2}$} & \multirow{2}{*}{$2$}\\
&\\
\multirow{2}{*}{39} & \multirow{2}{*}{} & \multirow{2}{*}{$A_2^{\oplus 3}$} & \multirow{2}{*}{$A_2$} & \multirow{2}{*}{$A_2^*\oplus\Z/3\Z$} & \multirow{2}{*}{$2$} & \multirow{2}{*}{$\frac{2}{3}$} & \multirow{2}{*}{$\frac{4}{3}$}\\
&\\
\multirow{2}{*}{40} & \multirow{2}{*}{} & \multirow{2}{*}{$A_2^{\oplus 2}\oplus A_1^{\oplus 2}$} & \multirow{2}{*}{$\langle 6\rangle^{\oplus 2}$} & \multirow{2}{*}{$\langle 1/6\rangle^{\oplus 2}$} & \multirow{2}{*}{$\frac{7}{3}$} & \multirow{2}{*}{$\frac{1}{2}$} & \multirow{2}{*}{$\frac{11}{6}$}\\
&\\
\multirow{2}{*}{41} & \multirow{2}{*}{} & \multirow{2}{*}{$A_2\oplus A_1^{\oplus 4}$} & \multirow{2}{*}{$\left(\begin{matrix}4 & -2\\-2 & 4\end{matrix}\right)$} & \multirow{2}{*}{$\frac{1}{6}\left(\begin{matrix}2 & 1\\1 & 2\end{matrix}\right)$} & \multirow{2}{*}{$\frac{8}{3}$} & \multirow{2}{*}{$\frac{1}{2}$} & \multirow{2}{*}{$\frac{13}{6}$}\\
&\\
\multirow{2}{*}{42} & \multirow{2}{*}{} & \multirow{2}{*}{$A_1^{\oplus 6}$} & \multirow{2}{*}{$A_1^{\oplus 2}$} & \multirow{2}{*}{${A_1^*}^{\oplus 2}\oplus(\Z/2\Z)^2$} & \multirow{2}{*}{$3$} & \multirow{2}{*}{$\frac{1}{2}$} & \multirow{2}{*}{$\frac{5}{2}$}\\
&\\
\hline
\multirow{2}{*}{43} & \multirow{2}{*}{$1$} & \multirow{2}{*}{$E_7$} & \multirow{2}{*}{$A_1$} & \multirow{2}{*}{$A_1^*$} & \multirow{2}{*}{$\frac{3}{2}$} & \multirow{2}{*}{$\frac{3}{2}$} & \multirow{2}{*}{$0$}\\
&\\
\multirow{2}{*}{44} & \multirow{2}{*}{} & \multirow{2}{*}{$A_7$} & \multirow{2}{*}{$A_1$} & \multirow{2}{*}{$A_1^*\oplus\Z/2\Z$} & \multirow{2}{*}{$2$} & \multirow{2}{*}{$\frac{7}{8}$} & \multirow{2}{*}{$\frac{11}{8}$}\\
&\\
\multirow{2}{*}{45} & \multirow{2}{*}{} & \multirow{2}{*}{$A_7$} & \multirow{2}{*}{$\langle 8\rangle$} & \multirow{2}{*}{$\langle 1/8\rangle$} & \multirow{2}{*}{$2$} & \multirow{2}{*}{$\frac{7}{8}$} & \multirow{2}{*}{$\frac{11}{8}$}\\
&\\
\multirow{2}{*}{46} & \multirow{2}{*}{} & \multirow{2}{*}{$D_7$} & \multirow{2}{*}{$\langle 4\rangle$} & \multirow{2}{*}{$\langle 1/4\rangle$} & \multirow{2}{*}{$\frac{7}{4}$} & \multirow{2}{*}{$1$} & \multirow{2}{*}{$\frac{3}{4}$}\\
&\\
\multirow{2}{*}{47} & \multirow{2}{*}{} & \multirow{2}{*}{$A_6\oplus A_1$} & \multirow{2}{*}{$\langle 14\rangle$} & \multirow{2}{*}{$\langle 1/14\rangle$} & \multirow{2}{*}{$\frac{31}{14}$} & \multirow{2}{*}{$\frac{1}{2}$} & \multirow{2}{*}{$\frac{12}{7}$}\\
&\\
\multirow{2}{*}{48} & \multirow{2}{*}{} & \multirow{2}{*}{$D_6\oplus A_1$} & \multirow{2}{*}{$A_1$} & \multirow{2}{*}{$A_1^*$} & \multirow{2}{*}{$2$} & \multirow{2}{*}{$\frac{3}{2}$} & \multirow{2}{*}{$\frac{1}{2}$}\\
&\\
\multirow{2}{*}{49} & \multirow{2}{*}{} & \multirow{2}{*}{$E_6\oplus A_1$} & \multirow{2}{*}{$\langle 6\rangle$} & \multirow{2}{*}{$\langle 1/6\rangle$} & \multirow{2}{*}{$\frac{11}{6}$} & \multirow{2}{*}{$\frac{1}{2}$} & \multirow{2}{*}{$\frac{4}{3}$}\\
&\\
\multirow{2}{*}{50} & \multirow{2}{*}{} & \multirow{2}{*}{$D_5\oplus A_2$} & \multirow{2}{*}{$\langle 12\rangle$} & \multirow{2}{*}{$\langle 1/12\rangle$} & \multirow{2}{*}{$\frac{23}{12}$} & \multirow{2}{*}{$\frac{2}{3}$} & \multirow{2}{*}{$\frac{5}{4}$}\\
&\\
\multirow{2}{*}{51} & \multirow{2}{*}{} & \multirow{2}{*}{$A_5\oplus A_2$} & \multirow{2}{*}{$A_1$} & \multirow{2}{*}{$A_1^*\oplus\Z/3\Z$} & \multirow{2}{*}{$\frac{13}{6}$} & \multirow{2}{*}{$\frac{2}{3}$} & \multirow{2}{*}{$\frac{3}{2}$}\\
&\\
\multirow{2}{*}{52} & \multirow{2}{*}{} & \multirow{2}{*}{$D_5\oplus A_1^{\oplus 2}$} & \multirow{2}{*}{$\langle 4\rangle$} & \multirow{2}{*}{$\langle 1/4\rangle\oplus\Z/2\Z$} & \multirow{2}{*}{$\frac{9}{4}$} & \multirow{2}{*}{$\frac{1}{2}$} & \multirow{2}{*}{$\frac{7}{4}$}\\
&\\
\multirow{2}{*}{53} & \multirow{2}{*}{} & \multirow{2}{*}{$A_5\oplus A_1^{\oplus 2}$} & \multirow{2}{*}{$\langle 6\rangle$} & \multirow{2}{*}{$\langle 1/6\rangle\oplus\Z/2\Z$} & \multirow{2}{*}{$\frac{5}{2}$} & \multirow{2}{*}{$\frac{1}{2}$} & \multirow{2}{*}{$2$}\\
&\\
\multirow{2}{*}{54} & \multirow{2}{*}{} & \multirow{2}{*}{$D_4\oplus A_3$} & \multirow{2}{*}{$\langle 4\rangle$} & \multirow{2}{*}{$\langle 1/4\rangle\oplus\Z/2\Z$} & \multirow{2}{*}{$2$} & \multirow{2}{*}{$\frac{3}{4}$} & \multirow{2}{*}{$\frac{5}{4}$}\\
&\\
\multirow{2}{*}{55} & \multirow{2}{*}{} & \multirow{2}{*}{$A_4\oplus A_3$} & \multirow{2}{*}{$\langle 20\rangle$} & \multirow{2}{*}{$\langle 1/20\rangle$} & \multirow{2}{*}{$\frac{11}{5}$} & \multirow{2}{*}{$\frac{3}{4}$} & \multirow{2}{*}{$\frac{29}{20}$}\\
&\\
\multirow{2}{*}{56} & \multirow{2}{*}{} & \multirow{2}{*}{$A_4\oplus A_2\oplus A_1$} & \multirow{2}{*}{$\langle 30\rangle$} & \multirow{2}{*}{$\langle 1/30\rangle$} & \multirow{2}{*}{$\frac{71}{30}$} & \multirow{2}{*}{$\frac{1}{2}$} & \multirow{2}{*}{$\frac{28}{15}$}\\
&\\
\multirow{2}{*}{57} & \multirow{2}{*}{} & \multirow{2}{*}{$D_4\oplus A_1^{\oplus 3}$} & \multirow{2}{*}{$A_1$} & \multirow{2}{*}{$A_1^*$} & \multirow{2}{*}{$\frac{5}{2}$} & \multirow{2}{*}{$\frac{1}{2}$} & \multirow{2}{*}{$2$}\\
&\\
\multirow{2}{*}{58} & \multirow{2}{*}{} & \multirow{2}{*}{$A_3^{\oplus 2}\oplus A_1$} & \multirow{2}{*}{$A_1$} & \multirow{2}{*}{$A_1^*\oplus\Z/4\Z$} & \multirow{2}{*}{$\frac{5}{2}$} & \multirow{2}{*}{$\frac{1}{2}$} & \multirow{2}{*}{$2$}\\
&\\
\multirow{2}{*}{59} & \multirow{2}{*}{} & \multirow{2}{*}{$A_3\oplus A_2\oplus A_1^{\oplus 2}$} & \multirow{2}{*}{$\langle 12\rangle$} & \multirow{2}{*}{$\langle 1/12\rangle\oplus\Z/2\Z$} & \multirow{2}{*}{$\frac{8}{3}$} & \multirow{2}{*}{$\frac{1}{2}$} & \multirow{2}{*}{$\frac{13}{6}$}\\
&\\
\multirow{2}{*}{60} & \multirow{2}{*}{} & \multirow{2}{*}{$A_3\oplus A_1^{\oplus 4}$} & \multirow{2}{*}{$\langle 4\rangle$} & \multirow{2}{*}{$\langle 1/4\rangle\oplus\Z/2\Z$} & \multirow{2}{*}{$3$} & \multirow{2}{*}{$\frac{1}{2}$} & \multirow{2}{*}{$\frac{5}{2}$}\\
&\\
\multirow{2}{*}{61} & \multirow{2}{*}{} & \multirow{2}{*}{$A_2^{\oplus 3}\oplus A_1$} & \multirow{2}{*}{$\langle 6\rangle$} & \multirow{2}{*}{$\langle 1/6\rangle\oplus\Z/3\Z$} & \multirow{2}{*}{$\frac{5}{2}$} & \multirow{2}{*}{$\frac{1}{2}$} & \multirow{2}{*}{$2$}\\
&\\
\hline
\caption{Mordell-Weil lattices of rational elliptic surfaces with Mordell-Weil rank $r\geq 1$.}
\end{longtable}
\end{center}

\bibliographystyle{alpha}
\bibliography{references_intersection_gaps}
\end{document}